\documentclass[11pt,letterpaper]{amsart}
\usepackage[utf8]{inputenc}
\usepackage[T1]{fontenc}
\usepackage[english]{babel}
\usepackage{subfigure}
\usepackage{amsmath}
\usepackage{amsthm} 
\usepackage{framed}
\usepackage[normalem]{ulem}
\usepackage{hyperref} 
\usepackage[OT2,OT1]{fontenc}
\usepackage{amsfonts}
\usepackage{tikz-cd}
\usepackage{tikz}
\usepackage{amssymb}
\usepackage{csquotes}
\usepackage{calrsfs} 
\usepackage[
backend=biber,
style=alphabetic,
sorting=nty,doi=false,isbn=false, url=false
]{biblatex}
\usepackage{accents}
\usepackage{adjustbox}
\usepackage{xcolor}
\definecolor{darkgreen}{RGB}{0,100,0}
\definecolor{DARKGREEN}{RGB}{0,100,0}

\definecolor{darkblue}{RGB}{0,0,139}
\definecolor{DARKBLUE}{RGB}{0,0,139}

\newlength{\dhatheight}

\makeatletter
\newcommand{\xleftrightarrow}[2][]{\ext@arrow 3359\leftrightarrowfill@{#1}{#2}}
\makeatother

\newcommand{\xdasharrow}[2][->]{
\tikz[baseline=-\the\dimexpr\fontdimen22\textfont2\relax]{
\node[anchor=south,font=\scriptsize, inner ysep=1.5pt,outer xsep=2.2pt](x){#2};
\draw[shorten <=3.4pt,shorten >=3.4pt,dashed,#1](x.south west)--(x.south east);
}
}
\theoremstyle{plain}
\newtheorem*{theorem*}{Theorem}
\newtheorem{theorem}{Theorem}[section]
\newtheorem{prop}[theorem]{Proposition}
\newtheorem{defi}[theorem]{Definition}
\newtheorem{lemma}[theorem]{Lemma}
\newtheorem{corollary}[theorem]{Corollary}
\theoremstyle{definition}
\newtheorem{example}[theorem]{Example}
\newtheorem{conj}[theorem]{Conjecture}

\newtheorem*{corollary*}{Corollary}
\newtheorem*{prop*}{Proposition}

\newtheoremstyle{named}{}{}{\itshape}{}{\bfseries}{.}{.5em}{\thmnote{#3's }#1}
\theoremstyle{named}

\newcommand{\C}{\mathbb{C}}
\newcommand{\Z}{\mathbb{Z}}
\newcommand{\Pro}{\mathbb{P}}
\newcommand{\Q}{\mathbb{Q}}
\newcommand{\A}{\mathbb{A}}
\newcommand{\N}{\mathbb{N}}

\newcommand{\et}{\textup{\'et}}
\newcommand{\ch}{\text{Chow}}

\newcommand{\alg}{\operatorname{alg}}

\newcommand{\CH}{\operatorname{CH}}
\newcommand{\tors}{\operatorname{tors}}
\theoremstyle{remark}
\newtheorem{remark}[theorem]{Remark}
\tikzcdset{
  cells={font=\everymath\expandafter{\the\everymath\displaystyle}},
}

\newcommand{\info}{{
  \bigskip
  \footnotesize
  \textsc{Departamento de Matemática, Universidad Técnica Federico Santa María, Avenida España 1680, Valparaíso, Chile}\par\nopagebreak
  \textit{E-mail address}: \texttt{pedro.montero@usm.cl }\par\nopagebreak
\vspace{2mm}
  \textsc{Université de Strasbourg, Institut de recherche mathématique avancée (IRMA), France}\par\nopagebreak
  \textit{E-mail address}: \texttt{ivan.rosas-soto@math.unistra.fr}
  }}
  

\author{Pedro Montero, Iv\'an Rosas-Soto}
\date{November 2026}
\title[Algebraic Cycles on Some Fano Varieties. Part II]{Algebraic cycles on some Fano varieties. Part~II: Hodge structures of level one}
\begin{document}
\maketitle

\begin{abstract}
We study Chow groups and étale motivic cohomology groups of smooth complete intersections with Hodge structures of level one, classified by Deligne and Rapoport, with particular attention to fivefolds. We extend these results to an étale motivic context and recover an analogous finite-dimensionality in the sense of Kimura.
\end{abstract}


\section*{Conventions and notations}
Throughout this paper we will assume that the base field $k$ is of $\text{char}(k)=0$, where in some subsections we will restrict to $k=\C$. For an abelian group $G$ we denote $G_{\tors}$ the subgroup of torsion elements of $G$. For a commutative ring $\A$ we denote $G_{\A}:=G\otimes_\Z \A$ and for an element $g \in G$ we denote $g_{\A}:= g\otimes 1 \in G_{\A}$. 

\section{Introduction}

Chow groups provide important information about the algebro-geometric structure of a variety, which can be used to address different problems in algebraic geometry, such as rationality questions, enumerative problems, and the behaviour of such properties in families, among others.

One of the open questions about algebraic cycles in a complex smooth projective variety $X$ is the so-called \emph{Hodge conjecture}, which relates the algebro-geometric information of the variety with its analytic counterpart. We refer the interested reader to the survey book \cite{Lew} for a comprehensive introduction to the subject. Originally, the version of the conjecture with integral coefficients was stated as follows: for all smooth projective varieties $X$ over $\C$ and for all $0 \leq i \leq \operatorname{dim}(X)$, the cycle class map
\[
    c^i:\CH^i(X) \to \operatorname{H}^{2i}_B(X,\Z)
\]
surjects onto the integral Hodge classes $\operatorname{Hdg}^{2i}(X,\Z)$. In this form, Atiyah and Hirzebruch \cite{AH} provided a counterexample in 1962, consisting of a torsion class in $H^4_B(X,\Z)$, for a Godeaux-Serre variety $X$ of dimension $7$, which is not algebraic. The situation does not improve if one considers only the torsion-free part of Hodge classes of a smooth projective variety: Kollár, in \cite{BCC}, produced the so-called \emph{Trento examples}, obtained by considering very general hypersurfaces of dimension $4$ and degree $125$, where the generator class of $\operatorname{H}^4_B(X,\Z)$ is not algebraic.

Despite this negative evidence for the integral version of the Hodge conjecture, its rational counterpart remains open to this day: for all smooth projective varieties $X$ over $\C$ and for all $0 \leq i \leq \operatorname{dim}(X)$ the cycle class map
\[
    c_\Q^i:\CH^i(X)_\Q \to \operatorname{Hdg}^{2i}(X,\Q)
\]
is surjective. There exists a generalised version, concerning Hodge structures of different weights, which may be stated as follows: let $X$ be a smooth projective variety, let $k \in \N$, and let $H \subset \operatorname{H}^k(X,\Q)$ be a sub-Hodge structure of coniveau $\geq c$ (see also \S \ref{section:Chow} and \cite[Lecture 7]{Lew}). Then there exists a closed algebraic subset $Y\subset X$ of codimension $c$ such that
\[
    H \subset \ker\big(j^*:\operatorname{H}_B^k(X,\Q) \to \operatorname{H}_B^k(X\setminus Y,\Q)\big).
\]

As we have mentioned, algebraic cycles and topological information also play a role in the study of rationality (or birational) problems. The first birational invariant of smooth projective varieties to recall is the Artin-Mumford invariant, given by $\operatorname{H}^3_B(X,\Z)_{\operatorname{tors}}$, which encodes the topological information of the Brauer group of $X$, in other words, it provides information about the Severi-Brauer varieties of $X$. Another invariant of topological nature is the obstruction to the integral Hodge conjecture in codimension $2$, i.e.\ $(\operatorname{Hdg}^4(X,\Z)/\operatorname{Im}(c^2))_{\operatorname{tors}}$. Similarly, for curve classes, one has that $\operatorname{Hdg}^{2n-2}(X,\Z)/\operatorname{Im}(c^{n-1})$ is a birational invariant of smooth projective varieties. All of these invariants are related to another broad family of birational invariants, the so-called \emph{unramified cohomology groups} introduced by Colliot-Thélène and Ojanguren in \cite{CTO89}, which have been successfully used to measure the failure of the Hodge conjecture for degree four integral classes \cite{CTV} and to produce new counterexamples to rationality questions \cite{CTP16} by means of the decomposition of the diagonal (see \S \ref{section:Chow}). We refer the interested reader to the survey articles \cite{Pey19,Voi19,Sch21} for an introduction to the subject.

\vspace{2mm}

In this article we study algebraic cycles on some Fano varieties whose middle cohomology has a particular Hodge structure, namely of level one, both classically and from the point of view of the étale topology, in the following sense.

\vspace{2mm}

We can obtain Chow groups by taking the hypercohomology of Bloch's complexes. More precisely, for each $j$, the presheaf $z^i(-,j)\colon U \mapsto z^i(U,j)$ associated with Bloch's complex is a sheaf for all topologies $t \in \left\{\et,\ \operatorname{Nis},\ \operatorname{Zar} \right\}$ (see \cite{Ge04}), and hence $z^i(-,\bullet)$ is a complex of sheaves on the small étale, Nisnevich, and Zariski sites of $X$. We define the complex of $t$-sheaves
\[
\Z_X(i)_t := z^i(-,\bullet)_t [-2i].
\]
In this case, by taking $t = \operatorname{Zar}$, the hypercohomology of the complex recovers Chow groups:
\[
\CH^i(X) \simeq \mathbb{H}^{2i}_{\operatorname{Zar}}(X, \Z(i)).
\]
From the perspective of triangulated categories, we can recover Chow groups by using the category of Nisnevich motives as
\[
\CH^i(X)\simeq \operatorname{Hom}_{\operatorname{DM}(k,\Z)}(M(X),\mathbf{1}(i)[2i]).
\]

All of these isomorphisms are related to the Zariski or Nisnevich topology. As we have mentioned, we may also proceed by considering the étale topology and obtain an étale version of Chow groups. Depending on the reference, this is called \emph{Lichtenbaum cohomology} or \emph{étale motivic cohomology}, which we denote by $\CH^i_\et(X)$. This new approach provides additional tools to study problems about algebraic cycles on smooth projective varieties, in particular the Hodge conjecture. For a smooth projective variety over a field of characteristic zero, we define
\[
    \CH^i_\et(X):=\operatorname{Hom}_{\operatorname{DM}_\et(k,\Z)}(M_\et(X),\Z(i)[2i])\simeq \mathbb{H}^{2i}_\et(X,\Z(i)).
\]

A first important feature (see \cite{Kahn,CTV}) of the étale Chow groups is that they fit into the following commutative diagram
\[
\begin{tikzcd}
    0 \arrow[r]& \CH^2(X) \arrow[r,"\kappa^2"]\arrow[d]& \CH^2_\et(X) \arrow[d] \arrow[r] & \operatorname{H}_{\operatorname{nr}}^3(X,\Q/\Z(2))\arrow[d] \arrow[r] & 0 \\
    0 \arrow[r] & \operatorname{Im}(c^2)\arrow[r] & \operatorname{H}^4_B(X,\Z)\arrow[r] & \operatorname{H}^4_B(X,\Z)/\operatorname{Im}(c^2)\arrow[r]& 0
\end{tikzcd}
\]
thus measuring the failure of the integral Hodge conjecture in degree $4$.

More generally, one of the fundamental results that illustrate the relation between algebraic cycles and étale motivic cohomology was given by Rosenschon and Srinivas in \cite{RS}, concerning the Hodge and Tate conjectures. In brief, for a smooth projective variety $X$ over $\C$, and for fixed $i \in \N$, there exists a cycle class map $c_\et^i:\CH_\et^i(X) \to \operatorname{H}_B^{2i}(X,\Z)$, whose image is contained in $\operatorname{Hdg}^{2i}(X,\Z)$. Therefore, we can formulate an étale analogue of the Hodge conjecture, namely: to ask whether the map $c_\et^i:\CH_\et^i(X) \to \operatorname{Hdg}^{2i}(X,\Z)$ is surjective. The first point to note is that, in contrast with the usual cycle class map, on torsion classes we have a surjection
\[
    c_\et^i\Big|_{\operatorname{tors}}:\CH_\et^i(X)_{\operatorname{tors}} \twoheadrightarrow \operatorname{Hdg}^{2i}(X,\Z)_{\operatorname{tors}}.
\]
This follows from the isomorphism between étale motivic cohomology groups with finite coefficients and étale cohomology, together with the comparison between étale and Betti cohomology with finite coefficients, see \cite{Mil80}. Consequently, we obtain the equivalence
\[
    c^i_\Q:\CH^i(X)_\Q \twoheadrightarrow \operatorname{Hdg}^{2i}(X,\Q) \iff c_\et^i:\CH_\et^i(X) \twoheadrightarrow \operatorname{Hdg}^{2i}(X,\Z).
\]

Since $c^i=c^i_\et\circ\kappa^i$, where $\kappa^i\colon\CH^i(X)\to\CH^i_\et(X)$ denotes the canonical comparison map associated with the change of sites (see \S\ref{section:Chow etale}), the integral Hodge conjecture for a variety satisfying the rational Hodge conjecture is thus controlled by the comparison maps $\kappa^i$, which are isomorphisms after tensoring with $\Q$. The étale point of view is also natural at the level of motives: a morphism of integral étale motives is an isomorphism if and only if it is so with rational coefficients and modulo every prime $\ell$, while $\operatorname{DM}_\et(k,\Z/\ell)$ is equivalent to the derived category of étale sheaves of $\Z/\ell$-modules and $\operatorname{DM}_\et(k,\Q)\simeq \operatorname{DM}(k,\Q)$ \cite{CD16}. This suggests that it may be easier to decompose integral étale motives than integral Chow motives, a point of view that we adopt in \S\ref{section:Motives}.

\vspace{2mm}

In connection with rationality problems, one of the major breakthroughs was obtained by Clemens and Griffiths \cite{CG72}, who proved that a smooth cubic threefold is unirational but not rational, relying heavily on the Hodge structure of these varieties (cf. the recent preprint \cite{EdGFS25}, where the authors prove that a very general cubic threefold is not even stably rational). More precisely, for such a threefold $X\subset \mathbb{P}^4$, the Hodge decomposition of its middle cohomology, $\operatorname{H}^3(X,\C) = V\oplus \overline{V}$, splits into two conjugate spaces. Consequently, as in the case of algebraic curves, one can associate to $X$ an \emph{intermediate Jacobian}, which is a principally polarised abelian variety of dimension $5$ in this case. The proof of Clemens and Griffiths proceeds by showing that this abelian variety is not the Jacobian of a curve. More generally, such a decomposition of the middle cohomology gives rise to a weight-one Hodge structure, which can be associated with a principally polarised abelian variety related to homologically trivial algebraic cycles, just as for an algebraic curve $C$, where the Jacobian coincides with the Albanese variety, the latter being dual to the group $\operatorname{Pic}^\circ(C)$ of degree-zero $0$-cycles.

On the other hand, from the perspective of $0$-cycles, an important class of algebraic varieties is given by Fano varieties (i.e.\ smooth projective varieties with ample anticanonical divisor). These form a concrete class of rationally connected varieties, meaning that all of their $0$-cycles are rationally equivalent, and they arise naturally in the birational classification of higher-dimensional algebraic varieties (see, for instance, \cite{Mur04} for a survey on algebraic cycles on Fano varieties). In regard of these considerations, in this article we focus on Fano varieties with Hodge structures of level one. Quartic double fivefolds, which are Fano fivefolds whose middle cohomology is numerically similar to that of a Calabi-Yau threefold, are studied in the first part of this work \cite{MRI}.

From a classification perspective, the simplest Fano varieties arise as hypersurfaces or, more generally, as smooth complete intersections of small degree in projective space. Regarding the classification of varieties of complete intersection with Hodge structures of level one given by Deligne and Rapoport in \cite{Del} and \cite{Rap}:
\begin{itemize}
    \item Fano threefolds given by:
    \begin{itemize}
        \item Cubic hypersurfaces in $\mathbb{P}^4$,
        \item Complete intersections of a quadric and a cubic in $\mathbb{P}^5$, and
        \item Quartic hypersurfaces in $\mathbb{P}^4$,
    \end{itemize}
    where the associated intermediate Jacobian is a principally polarised variety of dimension $5$, $20$ and $30$, respectively.
    \item Intersections of two quadrics in $\mathbb{P}^{2m+3}$. In this case, the intermediate Jacobian is the Jacobian of a hyperelliptic curve of genus $m+1$, see \cite{reid}.
    \item Intersections of three quadrics in $\mathbb{P}^{2m+4}$. In this case, the intermediate Jacobian is the Prym variety of a double cover of a plane curve of degree $2m+5$, see \cite{Beau}.
    \item If $X \subset \mathbb{P}^6$ is a smooth cubic hypersurface, then its intermediate Jacobian is a $21$-dimensional principally polarised abelian variety.
\end{itemize}

In this article, we push further the study of algebraic cycles on such varieties, which, \emph{grosso modo}, may be described as follows: for fivefolds, their Chow groups are fully characterised by their structure as complete intersections and their intermediate Jacobians. For further details about the classification of algebraic cycles of such varieties, we refer to \cite{reid, Beau, Mur04}, and more recently to the work of Fu and Tian \cite{FT}, who analysed the case of cubic fivefolds. This latter work provided the starting point for our study of algebraic cycles on Fano fivefolds from the point of view of classification of Fano varieties. For the threefolds and fivefolds in this list, we show that all the comparison maps $\kappa^i$ are isomorphisms (\S\ref{section:Chow etale level 1}), and that the middle part of their integral étale motive is, up to a Tate twist, the $h^1$ of the intermediate Jacobian (\S\ref{section:Motives}). Beyond complete intersections, higher-dimensional Fano manifolds of large index have been classified by Fujita and Mukai (see \cite{IP99}), and among them Gushel--Mukai varieties of odd dimension also carry Hodge structures of level one \cite{DK18,DK20} (see Remark~\ref{rmk:GM}).

The following result concerns the comparison maps $\kappa^i$ for the varieties with Hodge structure of level one of \S\ref{section:Chow etale level 1} and \S\ref{section:Motives}.

\begin{theorem}[{\S\ref{section:Chow etale level 1}, Remark~\ref{rmk:GM}}]\label{thm:main-kappa}
Let $X$ be a smooth projective variety over $\C$ which is a Fano threefold, a cubic fivefold, a fivefold complete intersection of two quadrics in $\mathbb{P}^7$ or of three quadrics in $\mathbb{P}^8$, or a Gushel--Mukai threefold or fivefold. Then the comparison maps $\kappa^i\colon \CH^i(X)\to\CH^i_\et(X)$ are isomorphisms for all $i\geq 0$.
\end{theorem}

\vspace{2mm}

The article is structured as follows. In \S \ref{section:Chow} and \S \ref{section:Chow etale} we recall the necessary preliminaries on Chow groups, intermediate Jacobians, decompositions of the diagonal and étale Chow groups. In \S \ref{section:Chow level 1} and \S \ref{section:Chow etale level 1} we describe the Chow and étale Chow groups of smooth complete intersections with Hodge structures of level one, with special emphasis on fivefolds, and in \S \ref{section:Motives} we show that some of these results lift to the motivic setting.

\section*{Acknowledgements}

The authors would like to thank the organisers of the conference ``Motivic Homotopy in Interaction'', held at CIRM Luminy from 4 to 8 November 2024. It was during this conference that this collaboration began, and we are especially grateful to \textsc{Fr\'ed\'eric D\'eglise}, \textsc{Adrien Dubouloz}, and \textsc{Johannes Nagel} for their constant support. We warmly thank \textsc{Stefan Schreieder} for all his comments and for pointing us to relevant references relating our results to his work on (refined) unramified cohomology groups. We are grateful to \textsc{Lie Fu} for helpful discussions. We also thank the anonymous referees for the careful reading of previous versions of this article and for helpful suggestions that improved several arguments and the exposition.

The first author was partially supported by ANID FONDECYT Regular grant 1231214 and by the ECOS-ANID Project ECOS230044. The second author was partially supported by the French ANR project “HQ-DIAG” (ANR-21-CE40-0015).

\section{Preliminaries on Chow groups}\label{section:Chow}

Let $X$ be a smooth projective variety over a field $k$ of dimension $n$. We denote the Chow group of $X$ in codimension $i$ by $\CH^i(X)$, i.e., the group of closed subvarieties of $X$ of codimension $i$ modulo rational equivalence. We can grade Chow groups of $X$ by dimension of subvarieties, in which case we have $\CH^{i}(X):= \CH_{n-i}(X)$. 

We freely use the standard properties of Chow groups of smooth projective varieties, for which we refer to \cite{Fulton}: pull-backs $f^\ast\colon\CH^i(Y)\to \CH^i(X)$ along arbitrary morphisms $f\colon X\to Y$, and push-forwards $f_\ast[Z]=\deg(Z/f(Z))[f(Z)]$ along proper ones \cite[\S 1.4 and Chapter 8]{Fulton}, the intersection product, which makes $\CH^\ast(X)=\bigoplus_{i=0}^n\CH^i(X)$ a commutative graded ring compatible with pull-backs, and the projection formula $f_\ast(Z\cdot f^\ast Z')=f_\ast(Z)\cdot Z'$ \cite[Chapter 8]{Fulton}, the localization exact sequence $\CH_i(Z)\to \CH_i(X)\to \CH_i(X\setminus Z)\to 0$ for a closed subset $Z\subset X$ \cite[\S 1.8]{Fulton}, and the action of a correspondence $Z\in\CH^t(X\times Y)$ on cycles $T\in \CH^i(X)$, given by $Z(T)=(\operatorname{pr}_Y)_\ast(Z\cdot(T\times Y))\in \CH^{i+t-n}(Y)$ with $n=\dim(X)$ \cite[\S 16.1]{Fulton}. Let us recall from \cite[\S 10.3]{Fulton}:

\begin{defi}
    We say that a cycle $Z \in \CH^i(X)$ is \textup{algebraically trivial} if there is a smooth irreducible curve $C$ and $W \in \CH^i(C\times X)$ and two points $a$, $b \in C$ such that $W(a)=0$ and $W(b)=Z$. We denote by $\CH^i(X)_{\alg}$ the subgroup of algebraically trivial cycles of codimension $i$ on $X$.
\end{defi}

When $k=\C$ there is a well-defined cycle map $c^i_{\A}:\CH^i(X)_\A \to \operatorname{H}^{2i}_B(X,\A)$ to Betti cohomology for $\A=\Z, \ \Q$, whose image is a subset of the Hodge classes
\begin{align*}
    \text{Hdg}^{2i}(X,\A):=\left\{ \alpha \in \operatorname{H}^{2i}(X,\A) \ | \  \alpha_\C \in \operatorname{H}^{i,i}(X)\right\}.
\end{align*}
We denote $\CH^i(X)_{\hom}:= \text{ker}(c^i)$, where $c^i:=c^i_\mathbb{Z}$. It is clear that $\CH^i(X)_{\text{alg}} \subset \CH^i(X)_{\text{hom}}$, and the quotient group is called the \textit{Griffiths group}, denoted by
\begin{align*}
    \text{Griff}^i(X):=\CH^i(X)_{\hom}/\CH^i(X)_{\alg}.
\end{align*}

\begin{defi}
Assume that $k=\bar{k}$, and let $\Omega/k$ be a field extension of infinite transcendence degree. A map 
$\CH^i(X_\Omega)_{\alg} \to T(\Omega)$ 
with $T$ a variety is said to be \emph{regular} if, for every smooth connected variety $S$ with base point $s_0 \in S(\Omega)$, and for every cycle $\Gamma \in \CH^i(S \times X)$, the composition
\begin{align*}
S(\Omega) &\to \CH^i(X_\Omega)_{\alg} \to T(\Omega) \\
s &\mapsto \Gamma((s) - (s_0))
\end{align*}
is a morphism of algebraic varieties.

We say that $\CH^i(X)_{\alg}$ is \emph{weakly representable} if there exists an abelian variety $A/k$ and a regular isomorphism of groups $\CH^i(X_{\Omega})_{\alg} \xrightarrow{\simeq} A(\Omega)$.
\end{defi}

Recall that a smooth projective variety $X$ of dimension $n$ over $k$ is \emph{rational} (resp.\ \emph{stably rational}) if $X$ (resp.\ $X\times\mathbb{P}^m_k$ for some $m\in\N$) is birational to projective space, \emph{unirational} if there is a dominant rational map $\mathbb{P}^N_k\dashrightarrow X$, and \emph{rationally connected} if two general points of $X$ can be connected by a rational curve. See e.g. \cite[Chapter 1]{KSC04} and \cite[Chapter IV]{Kol96}.

In the present article, we use properties of birational invariants related to algebraic cycles. The first and most important of these is that the group of zero-cycles of a smooth projective variety is a stably birational invariant. In particular, for a stably rational variety $X$, we have $\CH_0(X) = \mathbb{Z}$. The converse does not hold in general, but this condition can still be used to detect when a variety is not (stably) rational, through the so-called \emph{decomposition of the diagonal}. To that end, we introduce the notion of a universally trivial $\CH_0$-group.

\begin{defi}
Let $X$ be a smooth projective variety over $k$ of dimension $n$.
\begin{enumerate}
    \item We say that $X$ has a \textup{universally trivial $\CH_0$-group} if $X$ has a zero-cycle $x$ of degree 1 such that $\CH_0(X_K) = \mathbb{Z} x$ for any field extension $k \subset K$.
    
    \item We say that $X$ admits a \textup{Chow decomposition of the diagonal} if one can write
    \[
    \Delta_X = X \times x + Z \in \CH^n(X \times X),
    \]
    where $Z$ is a cycle on $X \times X$ supported on $D \times X$, with $D \subset X$ a proper closed algebraic subset.
\end{enumerate}
\end{defi}

These two conditions are equivalent (see eg. \cite[Proposition 3.11]{Voi19}). Zero-cycles also encode information about codimension 2 cycles and curves in $X$, depending on their support, as studied by Bloch and Srinivas in \cite{BS83}. To make this precise, let us consider a closed subvariety $Z \subset X$ such that for every field extension $L/k$, we have $\CH_0((X \setminus Z)_L) = 0$. Under these assumptions, by \cite[Proposition 1]{BS83}, there exist $N \in \mathbb{N}^{\geq 1}$, a divisor $D \subset X$, and cycles $\Gamma_1, \Gamma_2 \in \CH^n(X \times X)$, with $\Gamma_1$ supported on $X \times Z$ and $\Gamma_2$ supported on $D \times X$, such that\footnote{Note that if $N = 1$ and $Z = z$ is a point, we recover the classical decomposition of the diagonal.}
\[
N \cdot \Delta_X = \Gamma_1 + \Gamma_2 \in \CH^n(X \times X).
\]

Depending on the dimension of $Z$, we obtain the following consequences:

\begin{table}[h]
\centering
\begin{tabular}{|c|c|}
\hline
\textbf{dim $Z$} & \textbf{Conclusion about codimension 2 cycles on $X$} \\ \hline
$0$ & $\CH_0$-universally trivial \\ \hline
$\leq 1$ & $\CH^2(X)_{\text{alg}}$ is representable \\ \hline
$\leq 2$ & $\text{Griff}^2(X) = 0$, and there exists $N \in \mathbb{N}^{\geq 1}$ such that $N \cdot \text{Griff}_1(X) = 0$ \\ \hline
$\leq 3$ & The cycle class map $c^2: \CH^2(X) \twoheadrightarrow \text{Hdg}^4(X, \mathbb{Z})$ is surjective \\ \hline
\end{tabular}

\vspace{2mm}

\caption{Summary of results from \cite[Theorems 1 and 2]{BS83}.}
\end{table}

We will also use a generalization of the notion of decomposition of the diagonal, which involves cycles of higher dimension in $X$:

\begin{theorem}[{\cite[Theorem 3.18]{VoiB}}]\label{Teovoi}
Let $X$ be a connected smooth projective complex variety of dimension $n$, and let $c \in \mathbb{N}^{\geq 1}$ be such that for all $0 \leq k < c$, the cycle class map
\[
c^{n-k}_\Q : \CH_k(X)_\Q \hookrightarrow \operatorname{H}^{2(n-k)}(X, \Q)
\]
is injective. Then there exists $N \in \mathbb{N}^{\geq 1}$ such that
\[
N \cdot \Delta_X = Z_0 + Z_1 + \cdots + Z_{c-1} + Z' \in \CH^n(X \times X)_\Q,
\]
where each $Z_i$ is supported on $W_i \times W_i'$ with $\dim_\C(W_i) = i$ and $\dim_\C(W_i') = n - i$, and $Z'$ is supported on $X \times T$ with $T \subset X$ and $\textup{codim}_X(T) \geq c + 1$.
\end{theorem}

In the context described in the introduction, we will be interested in various Jacobians associated with Hodge structures. Let us first recall the notion of the level of a Hodge structure:

\begin{defi}
    Let $\displaystyle H_\mathbb{C}=\bigoplus_{p+q=k} H^{p,q}$ be a Hodge structure of weight $k$. The associated Hodge filtration is
    \[
        F^r(H)=\bigoplus_{p\geq r} H^{p,q}.
    \]
    The \emph{level} of a nonzero Hodge structure $H$ is the integer
    \[
        \ell := \max \{\, |p-q|,\;  H^{p,q}\neq 0\}.
    \]
    i.e., a Hodge structure $H$ of weight $k$ has level at most $k-2c$ if and only if $F^c H=H$.
\end{defi}

Recall that if $G$ is an abelian group and $X$ a smooth complex projective variety, the coniveau filtration on $\operatorname{H}^k(X,G)$ is defined by
\[
\operatorname{N}^c \operatorname{H}^k(X,G) = \sum_Z \ker \left\{\operatorname{H}^k(X,G) \to \operatorname{H}^k(X \setminus Z,G)\right\},
\]
where the sum ranges over all closed subvarieties $Z \subset X$ of codimension at least $c$ (see e.g. \cite{BO74}). With this notation, we introduce the various Jacobians that will play a role in our analysis.

\begin{defi}
    Let $H$ be a pure Hodge structure of odd weight $2m+1$ with $m \in \N$. The \emph{intermediate Jacobian} of $H$ is the complex torus
    \[
        J(H)=\overline{F^{m+1}H}/H_\Z.
    \]
    In particular, we introduce the following notations:
    \begin{enumerate}
        \item For $X \in \mathrm{SmProj}_\C$ and $H=\operatorname{H}^{2m+1}(X,\Q)$, the intermediate Jacobian of $X$ is
        \[
            J^{2m+1}(X) := J(H).
        \]
        \item For $H = \operatorname{N}^1\operatorname{H}^{2m+1}(X,\Q)$, we define the \emph{algebraic intermediate Jacobian} as
        \[
            J^{2m+1}_{\alg}(X) := J(H).
        \]
        \item Let $H^{2m+1}_{\mathrm{Hdg}}(X)$ be the maximal Hodge substructure of level $1$ contained in $\operatorname{H}^{2m+1}(X,\Q)$. Then we set $J^{2m+1}_H(X) := J\!\left(H^{2m+1}_{\mathrm{Hdg}}(X)\right)$.
        \item The \emph{Walker intermediate Jacobian} of $X$ is defined by
        \[
            J^{2m+1}_W(X) := J(\operatorname{N}^{m+1} \operatorname{H}^{2m+1}(X,\Q)).
        \]
    \end{enumerate}
\end{defi}

\begin{remark}
    \begin{enumerate}
        \item If $H = H^{2m+1}_{\mathrm{Hdg}}(X)$ or $H = \operatorname{N}^1\operatorname{H}^{2m+1}(X,\Q)$, then the Hodge--Riemann bilinear relations on $\operatorname{H}^{2m+1}(X,\Q)$ induce a polarization on $J(H)$, making it a principally polarized abelian variety. In particular, if $\operatorname{H}^{2m+1}(X,\Q)$ has level $1$, then $J^{2m+1}(X)$ is itself an abelian variety.
        \item There are natural inclusions of complex tori
        \[
            J_{\alg}^i(X)\hookrightarrow J^i_{H}(X) \hookrightarrow J^i(X).
        \]
        \item There exists an Abel--Jacobi map
        \[
            \Phi_X^i: \CH^i(X)_{\hom} \to J^{2i-1}(X).
        \]
        Restricting this map to $\CH^i(X)_{\alg}$ yields an isomorphism of abelian varieties
        \[
            \Phi^i_X(\CH^i(X)_{\alg}) \simeq J^{2i-1}_{\alg}(X).
        \]
        Moreover, the Abel--Jacobi map fits into the exact sequence
        \[
            0 \to J^{2i-1}(X) \xrightarrow{\Phi_X^i} \operatorname{H}^{2i}_D(X,\Z(i)) \to \mathrm{Hdg}^{2i}(X,\Z)\to 0,
        \]
        where $H^{2i}_D(X,\Z(i))$ denotes the Deligne cohomology group of $X$ (see \cite[Corollaire 12.27]{Voi}).
        \item If $C$ is a smooth projective curve, we omit the index and simply write $J(C)$ for its Jacobian.
    \end{enumerate}
\end{remark}

\section{Preliminaries on \'etale Chow groups}\label{section:Chow etale}

In this section, we recall the definition and fundamental properties of the étale analogues of Chow groups, and gather some general results concerning these groups for smooth complete intersections of Hodge level one.

\subsection{\'Etale Chow groups} We consider an étale version of Chow groups, namely the well-known Lichtenbaum cohomology groups, which are defined as the hypercohomology of the étale sheafification of Bloch’s cycle complex. These groups were characterized by Rosenschon and Srinivas in \cite{RS} using étale hypercoverings.

In this context, we let $\text{Sm}_k$ denote the category of smooth separated $k$-schemes over a field $k$. We denote by $z^i(X,\bullet)$ the cycle complex of abelian groups defined by Bloch:
\begin{align*}
z^i(X,\bullet): \cdots \to z^i(X,j) \to \cdots \to z^i(X,1)\to z^i(X,0) \to 0,
\end{align*}
where $z^i(X,j)$ is the subgroup generated by closed subvarieties of codimension $i$ of $X\times \Delta^j$ that intersect properly with $X \times F$ for each face $F$ of the standard simplex $\Delta^j$, and the differentials are given by the alternating sum of the pullbacks of the face maps.

We remind the reader that the homology groups of this complex define the higher Chow groups, namely $\text{CH}^i(X,2i-j) := \text{H}_j(z^i(X,\bullet))$, and in particular, $\CH^i(X) \simeq \CH^i(X,0)$.

Let us recall that $z^i(X,j)$ and the complex $z^i(X,\bullet)$ are covariantly functorial for proper morphisms and contravariantly functorial for flat morphisms between smooth $k$-schemes (see \cite[Proposition 1.3]{Blo}). Therefore, for any topology $t \in \left\{\text{fppf},\ \et,\ \text{Nis},\ \text{Zar} \right\}$, we obtain a complex of $t$-presheaves
\[
z^i(-,\bullet)\colon U \mapsto z^i(U,\bullet).
\]
In particular, for each $j$, the presheaf $z^i(-,j)\colon U \mapsto z^i(U,j)$ is a sheaf for all topologies $t \in \left\{\text{fppf},\ \et,\ \text{Nis},\ \text{Zar} \right\}$ (see \cite[Lemma 3.1]{Ge04}), and hence $z^i(-,\bullet)$ is a complex of sheaves on the small étale, Nisnevich, and Zariski sites of $X$.

We define the complex of $t$-sheaves
\[
R_X(i)_t := \left(z^i(-,\bullet)_t \otimes R \right)[-2i],
\]
where $R$ is an abelian group. In what follows, we will restrict to the cases $t = \text{Zar}$ or $t = \et$, and compute the corresponding hypercohomology groups $\mathbb{H}^j_t(X, R_X(i)_t)$.

For instance, taking $t = \text{Zar}$ and $R = \Z$, it follows from \cite[Proposition 4.2.9]{Voe00} that the hypercohomology of the complex recovers the higher Chow groups:
\[
\CH^i(X,2i-j) \simeq \mathbb{H}^j_{\text{Zar}}(X, \Z(i)).
\]
We denote the motivic and Lichtenbaum cohomology groups with coefficients in $R$ as
\[
\text{H}_M^j(X, R(i)) = \mathbb{H}^j_{\text{Zar}}(X, R(i)), \quad \text{H}_{M,\et}^j(X, R(i)) = \mathbb{H}^j_{\et}(X, R(i)),
\]
and in particular, we set $\CH^i_\et(X) := \text{H}^{2i}_{M,\et}(X, \Z(i))$.

Let $\pi \colon X_\et \to X_{\text{Zar}}$ be the canonical morphism of sites. The adjunction morphism $\Z_X(i) \to R\pi_* \pi^* \Z_X(i) = R\pi_* \Z_X(i)_\et$ induces the \emph{comparison morphisms}
\[
\kappa^{j,i} \colon \text{H}_M^j(X, \Z(i)) \longrightarrow \text{H}_{M,\et}^j(X, \Z(i))
\]
for all bidegrees $(j,i) \in \Z^2$.

Thanks to \cite[Theorem 6.18]{VV}, we know that the comparison map $\kappa^{j,i}$ is an isomorphism for $j \leq i+1$ and a monomorphism for $j \leq i+2$. For simplicity, we write 
\[
\kappa^i := \kappa^{2i,i}:\CH^i(X)\to \CH^i_\et(X)
\]
for the comparison map between Chow and Chow \'etale groups.

In some cases, it is possible to obtain more information about the Lichtenbaum cohomology groups and the comparison map between them and the higher Chow groups. For instance, there is a quasi-isomorphism $R_X(0)_\et = R$, where the right-hand side is regarded as an étale sheaf. Thus, the Lichtenbaum cohomology agrees with the usual étale cohomology, i.e.,
\[
\text{H}_{M,\et}^j(X, R(0)) \simeq \text{H}^j_\et(X, R)
\]
for all $j \in \mathbb{N}$, and in particular, $\text{CH}^0_\et(X) = \Z^{\pi_0(X)}$.

In the next case, $i = 1$, since there is a quasi-isomorphism of complexes $\Z_X(1)_\et \simeq \mathbb{G}_m[-1]$, we obtain the following isomorphisms:
\[
\text{CH}^1(X) \simeq \text{CH}_\et^1(X) = \text{Pic}(X), \quad
\text{H}_{M,\et}^3(X, \Z(1)) \simeq \text{H}_\et^3(X, \mathbb{G}_m[-1]) = \text{Br}(X),
\]
where $\text{Pic}(X)$ and $\text{Br}(X)$ denote the Picard and Grothendieck-Brauer groups of $X$, respectively. In fact, for bidegree $(j,1)$, by \cite[Corollary 3.4.3]{V}, there exists an isomorphism
\[
\text{H}_M^j(X, \Z(1)) \simeq \text{H}_{\text{Zar}}^{j-1}(X, \mathbb{G}_m),
\]
because the quasi-isomorphism $\Z_X(1) \simeq \mathbb{G}_m[-1]$ also holds in the Zariski topology.

As a particular case, consider
\[
\text{H}_M^3(X, \Z(1)) \simeq \text{H}_{\text{Zar}}^2(X, \mathbb{G}_m) = 0,
\]
since $\text{H}_M^j(X, \Z(i)) = 0$ for $j > 2i$, see \cite[Vanishing Theorem 19.3]{MVW}, whereas the Grothendieck-Brauer group $\text{Br}(X)$ is not always trivial (for instance, if $X$ is an Enriques surface, see e.g. \cite{Bea09}).

In bi-degree $(4,2)$, the comparison map is known to be injective but, in general, not surjective. More precisely, we have a short exact sequence
\[
0 \to \text{CH}^2(X) \xrightarrow{\kappa^2} \text{CH}^2_\et(X) \to \text{H}^3_{\text{nr}}(X, \Q/\Z(2)) \to 0,
\]
where $\mathcal{H}^3_\et(\Q/\Z(2))$ is the Zariski sheaf associated to the presheaf $U \mapsto \text{H}_\et^3(U, \Q/\Z(2))$, and the unramified part is given by the global sections $\text{H}_{\text{nr}}^3(X, \Q/\Z(2)) = \Gamma(X, \mathcal{H}^3_\et(\Q/\Z(2)))$. For a proof, see \cite[Proposition 2.9]{Kahn}.

If $k = \C$, the latter group surjects onto the torsion subgroup of the obstruction, in codimension 2, to the integral Hodge conjecture, i.e.,
\[
\text{H}_{\text{nr}}^3(X, \Q/\Z(2)) \twoheadrightarrow 
\left( \text{Hdg}^4(X, \Z) / \operatorname{Im}\left\{ c^2: \CH^2(X) \to \text{H}_B^4(X, \Z(2)) \right\} \right)_{\text{tors}},
\]
and hence, in general, it is nonzero and the comparison map $\kappa^2$ is not surjective. For further details, see \cite[Théorème 3.7]{CTV}.

\begin{remark}
For rational coefficients, the adjunction morphism
\[
\Q_X(i) \to R\pi_* \Q_X(i)_\et
\]
turns out to be an isomorphism (see \cite[Théorème 2.6]{Kahn}). Therefore, we have
\[
\text{H}_M^j(X, \Q(i)) \simeq \text{H}_{M,\et}^j(X, \Q(i))
\]
for all $(j,i) \in \Z^2$. In particular,
\[
\text{CH}^i(X)_\Q \simeq \text{CH}_\et^i(X)_\Q
\]
for all $i \in \N$.
\end{remark}

If $R$ is a torsion group, then the Lichtenbaum cohomology can be computed via étale cohomology. More precisely, let $\ell$ be a prime number, $r \in \N^{\geq 1}$, and $R = \Z/\ell^r$. Then there is a quasi-isomorphism
\[
(\Z/\ell^r)_X(i)_\et \xrightarrow{\sim} \mu_{\ell^r}^{\otimes i}.
\]

Passing to the direct limit, we also obtain a quasi-isomorphism
\[
(\Q_\ell/\Z_\ell)_X(i)_\et \xrightarrow{\sim} \varinjlim_r \mu_{\ell^r}^{\otimes i},
\]
and finally we define
\[
(\Q/\Z)_X(i)_\et := \bigoplus_\ell (\Q_\ell/\Z_\ell)_X(i)_\et \xrightarrow{\sim} \Q/\Z(i)_\et.
\]

The following result is a well-known consequence, referred to as the \emph{Suslin rigidity theorem}, concerning the morphism $\Z_X(i) \to R\pi_* \Z_X(i)_\et$ for $i \geq \dim(X)$ over $k = \bar{k}$.

\begin{prop}[{\cite[Theorem 4.2]{sus}, \cite[Section 2]{Geis}}]\label{lemGe}
Let $X$ be a smooth projective variety of dimension $n$ over an algebraically closed field $k$. Then, for $i \geq n$, the canonical morphism of sites $\pi \colon X_\et \to X_{\text{Zar}}$ induces a quasi-isomorphism of complexes of Zariski sheaves:
\[
\Z_X(i) \to R\pi_* \Z_X(i)_\et.
\]
In particular, $\operatorname{CH}^n(X) \simeq \operatorname{CH}_\et^n(X)$.
\end{prop}

Throughout the present article, we will use various cycle class maps, ranging from étale motivic and Lichtenbaum cohomology groups to Betti and étale cohomology:
\[
\begin{aligned}
&c^i\colon \CH^i(X) \to \text{H}^{2i}_B(X,\Z(i)), \quad (\text{resp. } c^i_{\Q}\colon \CH^i(X)_{\Q} \to \text{H}^{2i}_B(X,\Q(i))) \\
&c^i_{\et}\colon \CH^i_\et(X) \to \text{H}^{2i}_B(X,\Z(i)),
\end{aligned}
\]
the cycle class map and the Lichtenbaum cycle class map (see \cite{RS}). These maps satisfy the property that, as in the classical case of Chow groups, we have
\[
\operatorname{Im}(c^i_\et) \subset \operatorname{Hdg}^{2i}(X, \Z),
\quad \text{but} \quad
\operatorname{Im}(c^i_\et)_{\text{tors}} \simeq \operatorname{Hdg}^{2i}(X, \Z)_{\text{tors}}.
\]

\begin{remark}
Let us briefly discuss the functoriality of étale Chow groups: Since we work over a base field $k$ with $\operatorname{char}(k) = 0$, we can use an equivalent definition (see \cite[Theorem 7.1.2]{CD16}) of étale Chow groups in terms of the triangulated category of étale motives $\operatorname{DM}_\et(k, \Z)$. For more details, we refer the reader to \cite{Ayo} and \cite{CD16}. One defines the \textit{étale motivic cohomology} group of bidegree $(j,i)$ as
\[
\text{H}_{M,\et}^j(X, \Z(i)) := \operatorname{Hom}_{\operatorname{DM}_\et(k, \Z)}(M_\et(X), \Z(i)[j]).
\]
Thanks to the six functor formalism (see \cite{CD19}), we can recover the standard formalism for Chow groups (pull-backs, push-forwards of proper maps, intersection products, and the projection formula), now in the context of étale Chow groups. See \cite{RoSo} for further details on the construction and compatibility properties.
\end{remark}

An \emph{étale correspondence} from $X$ to $Y$ is a cycle $Z \in \CH_\et(X \times Y)$. A correspondence $Z \in \CH_\et^t(X \times Y)$ acts on cycles in $X$ via
\[
Z(T) := (\operatorname{pr}_Y)_*(Z \cdot (T \times Y)) \in \CH_\et^{i + t - n}(Y),
\]
where $T \in \CH_\et^i(X)$ and $n = \dim(X)$.

\begin{defi}
We say that a cycle $Z \in \CH_\et^i(X)$ is \emph{algebraically trivial} if there exists a smooth irreducible curve $C$, a cycle $W \in \CH_\et^i(C \times X)$, and two points $a, b \in C$ such that $W(a) = 0$ and $W(b) = Z$. We denote the subgroup of algebraically trivial cycles of codimension $i$ by $\CH_\et^i(X)_{\alg}$.
\end{defi}

We denote by $\CH_\et^i(X)_{\text{hom}} := \ker(c^i_\et)$. It is clear that
\[
\CH_\et^i(X)_{\text{alg}} \subset \CH_\et^i(X)_{\text{hom}},
\]
although it is not always obvious that this inclusion is strict. However, $\CH_\et^i(X)_{\text{alg}}$ is always a divisible group, since $\CH_\et^i(X)_{\text{hom}}$ is the maximal divisible subgroup of $\CH_\et^i(X)$ \cite[\S 3]{Gei17}. The corresponding quotient is called the \textit{étale Griffiths group}, and is denoted as
\[
\operatorname{Griff}_\et^i(X) := \CH_\et^i(X)_{\text{hom}} / \CH_\et^i(X)_{\text{alg}}.
\]

\subsection{\'Etale motivic cohomology of complete intersections} Let us consider $X = X_{(d_1, \ldots, d_r)} \subseteq \mathbb{P}^{n+r}$ to be a smooth complete intersection of multidegree $(d_1, \ldots, d_r)$ and dimension $n$, defined over $\mathbb{C}$.

Since the cohomology of $X$ is torsion-free and coincides with that of $\mathbb{P}^n$ outside the middle degree by the Lefschetz hyperplane theorem, the only non-trivial Betti number of $X$ is $b_n(X)$, which is determined by the Euler characteristic $\chi(X)$ and thus can be computed from Whitney's formula for the total Chern class of the tangent bundle of $X$. Explicit formulas for $\chi(X)$ and $b_n(X)$ in terms of $n$ and $d_1,\ldots,d_r$ can be found, for instance, in \cite[Chapter 5]{Dim92}. For example, $b_3(X)=10$ for cubic threefolds and $b_5(X)=42$ for cubic fivefolds.

\begin{prop}
Let $X = X_{(d_1, \ldots, d_r)} \subset \mathbb{P}^{n+r}$ be a smooth complete intersection. Then all the groups $\CH^i_\et(X)$ are torsion-free, except possibly when $n = 2m + 1$ is odd and $i = m + 1$.
\end{prop}

\begin{proof}
Assume we are not in the exceptional case $\dim(X) = 2m + 1$ and $i = m + 1$. The cohomology groups of $X$ are torsion-free in every degree by Lefschetz hyperplane section theorem, hence so are their integral Hodge structures. It follows that
\[
\CH^i_\et(X)_{\text{tors}} \simeq \CH^i_\et(X)^{\text{hom}}_{\text{tors}},
\]
independently of $i$. In particular, we have an isomorphism
\[
\CH^i_\et(X)_{\text{tors}} \simeq J^{2i-1}(X)_{\text{tors}},
\]
where $J^i(X)$ is the intermediate Jacobian. By our assumptions on $i\neq m+1$, we have $\text{H}^{2i+1}(X, \Z) = 0$, and thus $J^{2i+1}(X) = 0$. Hence, $\CH^i_\et(X)^{\text{hom}}$ is torsion-free for every $i\neq m+1$. Since this group is divisible\footnote{In fact, it is the maximal divisible subgroup of $\CH^i_\et(X)$.}, it must be uniquely divisible, i.e., a $\Q$-vector space. Therefore, the short exact sequence
\[
0 \to \CH^i_\et(X)_{\text{hom}} \to \CH^i_\et(X) \to \operatorname{Im}(c^i_\et) \to 0
\]
splits, and using the isomorphism $\CH^i_\et(X)_\Q \simeq \CH^i(X)_\Q$, we conclude:
\[
\CH^i_\et(X) \simeq \operatorname{Im}(c^i_\et) \oplus \CH^i_{\text{hom}}(X)_\Q,
\]
where $\operatorname{Im}(c^i_\et)_{\text{tors}} \simeq \operatorname{Hdg}^{2i}(X, \Z)_{\text{tors}}=0$ in this case. 
\end{proof}

\begin{remark}
\begin{enumerate}
    \item In the excluded case $\dim(X) = 2m+1$ and $i = m+1$, we have:
    \[
    \CH_\et^{m+1}(X)_{\text{tors}} \simeq J^{2m+1}(X)_{\text{tors}} \simeq (\Q/\Z)^{\oplus b_{2m+1}(X)}.
    \]
    
    \item If $\dim(X) = 2m+1$, then the Hodge conjecture holds (rationally). By the equivalence in \cite{RS}, we have:
    \[
    0 \to \CH^i_\et(X)_{\text{hom}} \to \CH^i_\et(X) \to \Z[h^i] \to 0,
    \]
    which splits via a natural section. The same holds in even dimension, except possibly in the middle degree, where the Hodge conjecture remains open.
\end{enumerate}
\end{remark}

\section{Chow groups of complete intersections with Hodge level one}\label{section:Chow level 1}

Smooth complete intersections in the complex projective space whose Hodge structures have level one were fully classified by Deligne and Rapoport in \cite{Del} and \cite{Rap}, respectively. These varieties are the following:
\begin{enumerate}
    \item Fano threefolds given by:
    \begin{itemize}
        \item Cubic hypersurfaces in $\mathbb{P}^4$,
        \item Quartic hypersurfaces in $\mathbb{P}^4$, and
        \item Complete intersections of a quadric and a cubic in $\mathbb{P}^5$.
    \end{itemize}
    
    \item Intersections of two quadrics in $\mathbb{P}^{2m+3}$.
    
    \item Intersections of three quadrics in $\mathbb{P}^{2m+4}$.
    
    \item Cubic fivefolds in $\mathbb{P}^6$.
\end{enumerate}

Recall that if $G$ is an abelian group and $X$ a complex algebraic variety, then $\mathcal{H}^q(G)$ denotes the Zariski sheaf associated to the presheaf $U \mapsto \operatorname{H}_B^q(U,G)$. The unramified cohomology groups are defined as
\[
\operatorname{H}_{\operatorname{nr}}^i(X,G) := \operatorname{H}_{\operatorname{Zar}}^0(X,\mathcal{H}^i(G)),
\]
which are birational invariants and relate to the torsion in Chow groups. More generally, if $\mathcal{G}$ is an étale sheaf on $X$, we set
\[
\operatorname{H}_{\operatorname{nr}}^i(X,\mathcal{G}) := \operatorname{H}^0_{\operatorname{Zar}}(X,\mathcal{H}^i(\mathcal{G})),
\]
where $\mathcal{H}^i(\mathcal{G})$ is the Zariski sheaf associated to $U \mapsto \operatorname{H}^i_{\et}(U,\mathcal{G})$.

\subsection{Chow groups of an intersection of two quadrics}\label{subsec:Chow 2 quadrics}

Let $X = X_{(2,2)} \subset \mathbb{P}^{2m+3}$. When $m \geq 1$, the variety $X$ is rational: by \cite[Theorem 2.1]{DM98} we know that $X$ contains a line $\ell \simeq \mathbb{P}^1 \subset X$, and more generally, $X$ contains $j$-planes $H\simeq \mathbb{P}^j$ for all $1 \leq j \leq m$. Here, the linear projection $\pi_\ell:X \dashrightarrow \mathbb{P}^{2m+1}$ from a line induces a birational map (cf. \cite{HT21a}).

\begin{theorem}\label{tors}
Let $X = X_{(2,2)} \subset \mathbb{P}^{2m+3}$ with $m \geq 1$. Fix an integer $p \neq m$. Then for all $k \geq 1$, the groups $\operatorname{H}^p(X, \mathcal{H}^{p+k}(\Z))$ are torsion.
\end{theorem}

\begin{proof}
The result holds for $p = 0$ and $p = 1$. Let us now fix $p \leq m - 1$. The hypotheses of Theorem~\ref{Teovoi} are satisfied for $X$ with $c = p + 1$, so we have a decomposition:
\[
N_p \Delta_X = Z_0 + Z_1 + \cdots + Z_p + Z',
\]
where $Z'$ is supported on $X \times T$ with $\operatorname{codim}_X(T) \geq p + 1$. Denote by $\widetilde{W}_i$ and $\widetilde{T}$ the desingularizations of $W_i$ and $T$, respectively.

Since the cycle $\Delta_X$ acts as the identity on cohomology, we have that $N_p \Delta_X^\ast$ acts on $\text{H}^p(X, \mathcal{H}^q(\Z))$ as $N_p \cdot \operatorname{Id}$. Thanks to the decomposition of the diagonal, the action of the diagonal (times the constant $N_p$) factors through the following terms:
\begin{itemize}
    \item $\text{H}^p(\widetilde{W}_i, \mathcal{H}^{p+k}_{\widetilde{W}_i}(\Z))$, which vanishes since $p+k > \dim(\widetilde{W}_i) = i$, so $\mathcal{H}^{p+k}_{\widetilde{W}_i}(\Z) = 0$.
    
    \item $\text{H}^{p-c}(\widetilde{T}, \mathcal{H}^{p+k - c}_{\widetilde{T}}(\Z))$, where $c = \operatorname{codim}_X(T)$. As long as $p - c < 0$ we have that the latter group vanishes.
\end{itemize}

Therefore, $N_p \cdot \Delta_X$ acts as zero on $\text{H}^p(X, \mathcal{H}^{p+k}(\Z))$, and the result follows.

Now suppose $p \geq m + 1$. Set $c_p := 2m + 1 - p$ so that the hypotheses are verified for $c = c_p + 1$. Then, the $N_{c_p}\Delta_X^\ast$ factors through:
\begin{itemize}
    \item $\text{H}^{p - d_i}(\widetilde{W}_i, \mathcal{H}^{q - d_i}_{\widetilde{W}_i}(\Z))$, where $d_i = 2m + 1 - i$ and again $p - d_i < 0$, so these terms vanish.
    
    \item $\text{H}^p(\widetilde{T}, \mathcal{H}^{p+k}_{\widetilde{T}}(\Z))$, and since $\operatorname{codim}_X(T) \geq c_p + 1$, we get $\dim(T) \leq p - 1$, so again the cohomology group vanishes.
\end{itemize}
As before, we deduce that the groups $\text{H}^p(X, \mathcal{H}^{p+k}(\Z))$ are torsion.
\end{proof}

\begin{corollary}\label{zero}
Let $X = X_{(2,2)} \subset \mathbb{P}^{2m+3}$ with $m \geq 1$. Then, for any fixed integer $p \neq m$ and for all $k \geq 1$, we have
\[
\operatorname{H}^p(X, \mathcal{H}^{p+k}(\Z)) = 0.
\]
\end{corollary}

\begin{proof}
Let $\phi: X \dashrightarrow \mathbb{P}^{2m+1}$ be a birational map which is an isomorphism onto its image outside a line $\ell \subset X$, where $\operatorname{codim}_X(\ell) = 2m$ (i.e., $\dim(\ell) = 1$). Then:
\[
\Delta_X = \phi^*(\Delta_{\mathbb{P}^{2m+1}}) + Z \in \CH^{2m+1}(X \times X),
\]
where $Z$ is a cycle supported on $H\times H$, since $\phi$ is an isomorphism outside a linear subspace $H\subset X$ of $\operatorname{codim}_X(H)=m+1$. Indeed, the decomposition of the diagonal of $\Delta_{\mathbb{P}^{2m+1}}$ gives $N_p = 1$. Finally, the result follows directly from Theorem~\ref{tors}.
\end{proof}

\begin{corollary}
    Let $X=X_{(2,2)}\subset \mathbb{P}^{2m+3}$ with $m\geq 1$. Then
    \begin{align*}
        \CH^i(X)\simeq \left\{ \begin{array}{cl}
           \mathbb{Z}  & \textup{if }i\neq m+1, \\
            \mathbb{Z}\oplus J^{2m+1}(X) & \textup{if }i=m+1. 
        \end{array} \right.
    \end{align*}
\end{corollary}
\begin{proof}
    Let $n\in\N^{\geq 2}$. From the short exact sequence
    \begin{align*}
        0\to \operatorname{H}^{i-1}(X, \mathcal{H}^{i}(\Z))/n\to \operatorname{H}^{i-1}(X, \mathcal{H}^{i}(\Z/n)) \to \operatorname{H}^i(X, \mathcal{H}^{i}(\Z))[n] \to 0
    \end{align*}
and the fact that $\CH^i(X)/\textup{alg}\simeq  \operatorname{H}^i(X, \mathcal{H}^{i}(\Z))$ are torsion free, thus we have an isomorphism 
\begin{align*}
    \operatorname{H}^{i-1}(X, \mathcal{H}^{i}(\Z))/n\xrightarrow{\simeq}\operatorname{H}^{i-1}(X, \mathcal{H}^{i}(\Z/n)).
\end{align*}
By \cite[Proposition 1]{CSS}, there exists a surjection of abelian groups
\begin{align*}
    \operatorname{H}^{i-1}(X, \mathcal{H}^{i}(\Z/n)) \twoheadrightarrow \CH^i(X)[n]
\end{align*}
By Corollary \ref{zero}, for all $p\neq m$ and $k\geq 1$ the group $\operatorname{H}^p(X, \mathcal{H}^{p+k}(\Z)) = 0$, therefore for all $i\neq m+1$ one gets $\CH^i(X)[n]=0$ and consequently they are torsion free. The validity of the integral Hodge conjecture and the isomorphisms (see e.g. \cite[Theorem 1.2]{BT16})
    \begin{align*}
        \CH^i(X)_\Q\simeq \left\{ \begin{array}{cl}
           \mathbb{Q}  & \textup{if }i\neq m+1, \\
            \mathbb{Q}\oplus J^{2m+1}(X) & \textup{if }i=m+1. 
        \end{array} \right.
    \end{align*}
imply the result.
\end{proof}

\subsection{Chow groups of an intersection of three quadrics}\label{subsec:Chow 3 quadrics}

If we consider $X_{(2,2,2)} \subset \mathbb{P}^{2m+4}$, then for $m \geq 2$, the variety $X$ is rational (cf. \cite{HT21}). Moreover, invoking \cite[Theorem 2.1]{DM98}, we can observe that it contains a linear space $H \simeq \mathbb{P}^i \hookrightarrow X$, with $i = m$ if $m \leq 2$ and $i = m - 1$ if $m \geq 3$.

\subsubsection{Case when \texorpdfstring{$m = 2$}{}}

Let $X = X_{(2,2,2)} \subset \mathbb{P}^8$ be a smooth complete intersection of dimension 5. It contains a plane $H \simeq \mathbb{P}^2$ and admits a birational map $X \dashrightarrow \mathbb{P}^5$. We have the following information about its Chow groups:
\begin{itemize}
    \item $\CH^0(X) \simeq \Z$.
    \item $\CH^5(X) \simeq \Z$ since $X$ is rational. In particular, this implies that $X$ is $\CH_0$-universally trivial.
    \item $\CH^1(X) \simeq \Z[h]$, generated by a hyperplane section $h = c_1(\mathcal{O}_X(1))$.
    \item $\CH^2(X) \simeq \Z[h^2]$. Indeed, the integral Hodge conjecture holds for $X$ since it is $\CH_0$-universally trivial (see \cite[Theorem 1]{BS83}). Moreover, $\CH^2(X)_{\text{hom}} \simeq \CH^2(X)_{\text{alg}}$ and $\CH^2(X)_{\text{alg}} = 0$, hence the conclusion follows.
\end{itemize}

Regarding the group $\CH^4(X)$ of 1-cycles, we note: by \cite[Theorem 2]{BS83}, $\operatorname{Griff}^4(X)$ is torsion and annihilated by the integer $N$ in the decomposition of the diagonal. Since $X$ is $\CH_0$-universally trivial, $N = 1$, so $\CH^4(X)_{\text{alg}} = \CH^4(X)_{\text{hom}}$. Moreover, $\CH^4(X)_{\text{hom}}$ is torsion. We have a diagram
\[
\begin{tikzcd}[column sep=small]
0 \arrow{r} & \CH^4(X)_{\text{tors}} \arrow{r}\arrow{d} & \CH^4(X) \arrow{r}\arrow{d} & \CH^4(X)_\Q \arrow{r}\arrow{d} & \CH^4(X)_{\Q/\Z} \arrow{d}{\simeq} \arrow{r} & 0 \\
0 \arrow{r} & 0 \arrow{r} & \text{H}^{8}_B(X,\Z) \arrow{r} &  \text{H}^{8}_B(X,\Q) \arrow{r} &  \text{H}^{8}_B(X,\Q/\Z ) \arrow{r} & 0
\end{tikzcd}
\]
since $\CH^4(X)_{\Q/\Z} \simeq A^4(X)_{\Q/\Z} \simeq \Q/\Z$ and $\CH^4(X)_{\text{alg}} = \CH^4(X)_{\text{hom}}$, where $A^i(X) = \CH^i(X)/\CH^i(X)_{\text{alg}}$. The upper exact sequence implies, by analyzing the rank of $\CH^4(X)_\Q$, that $\CH^4(X) \simeq \Z \oplus \CH^4(X)_{\text{tors}}$ and then $\CH^4(X)_{\text{tors}} \simeq \CH^4(X)_{\text{hom}}$. Now, by \cite[Proposition 1]{CSS}, the map
\[
\text{H}^3_{\text{Zar}}(X,\mathcal{H}^4(\mu_{\ell^r}^{\otimes 4})) \to \CH^4(X)[\ell^r]
\]
is surjective. Using the Bloch-Ogus spectral sequence \cite{BO74},
\[
E^{p,q}_2 = \text{H}^p_{\text{Zar}}(X,\mathcal{H}^q(\mu_{\ell^r}^{\otimes 4})) \Longrightarrow \text{H}^{p+q}_\et(X,\mu_{\ell^r}^{\otimes 4}),
\]
and since $\text{H}^7_\et(X,\mu_{\ell^r}^{\otimes 4}) = 0$, all the graded pieces associated to $E^7_\infty$ vanish. We also have that
\[
E_\infty^{3,4} \simeq E_3^{3,4} \simeq \frac{E_2^{3,4}}{\operatorname{Im}(E_2^{1,5} \to E_2^{3,4})}.
\]
Moreover, we have a short exact sequence:
\[
0 \to \text{H}^5_{\text{nr}}(X,\Z(4))/\ell^r \to \text{H}^5_{\text{nr}}(X,\mu_{\ell^r}^{\otimes 4}) \to \text{H}^1_{\text{Zar}}(X,\mathcal{H}^5(\Z(4)))[\ell^r] \to 0.
\]
Due to the rationality of $X$, we have that $\text{H}^1_{\text{Zar}}(X,\mathcal{H}^5(\Z(4)))_{\text{tors}} = 0$, so 
\[
\text{H}^1_{\text{Zar}}(X,\mathcal{H}^5(\Z(4)))/\ell^r \simeq \text{H}^1_{\text{Zar}}(X,\mathcal{H}^5(\mu_{\ell^r}^{\otimes 4})).
\]
Using the same arguments as in \cite[Lemme 21]{FT}, one can deduce directly that $\text{H}^1_{\text{Zar}}(X,\mathcal{H}^5(\Z(4))) = 0$, and hence $E_2^{1,5} = 0$. It follows that $E_\infty^{3,4} = E_2^{3,4}$, and thus
\[
E_\infty^{3,4} = \text{H}^3_{\text{Zar}}(X, \mathcal{H}^4(\mu_{\ell^r}^{\otimes 4})),
\]
which implies $\CH^4(X)[\ell^r] = 0$.

\begin{prop}
Let $X = X_{(2,2,2)} \subset \mathbb{P}^8$ be a smooth complete intersection. Then:
\[
\operatorname{N}^2\operatorname{H}^5(X,\Z) := \sum_{\operatorname{codim}_X(X \setminus U) \geq 2} \ker \left\{ \operatorname{H}^5(X,\Z) \xrightarrow{j^*} \operatorname{H}^5(U,\Z) \right\} \simeq \operatorname{H}^5(X,\Z).
\]
\end{prop}

\begin{proof}
The variety $X$ contains a plane $H \simeq \mathbb{P}^2$ by \cite[Theorem 2.1]{DM98}, inducing a birational map $\phi : X \dashrightarrow \mathbb{P}^5$ (cf. \cite{HT21}). Let $\Gamma_\phi \subset X \times \mathbb{P}^5$ and $\Gamma_{\phi^{-1}} \subset \mathbb{P}^5 \times X$ be the closures of the graphs of $\phi$ and $\phi^{-1}$.

Since we are considering smooth projective varieties, $\phi$ is defined on an open $U \subset X$ with $\operatorname{codim}_X(X \setminus U) \geq 2$. Similarly, $\phi^{-1}$ is defined on $V \subset \mathbb{P}^5$ with $\operatorname{codim}_{\mathbb{P}^5}(\mathbb{P}^5 \setminus V) \geq 2$.

Let $Z = X \setminus U$ and $W = \mathbb{P}^5 \setminus V$. Then (see \cite[Lemme 3.5]{CTV}),
\[
\Gamma_{\phi^{-1}} \circ \Gamma_\phi = \Delta_X + \mathcal{Z}, \quad \Gamma_\phi \circ \Gamma_{\phi^{-1}} = \Delta_{\mathbb{P}^5} + \mathcal{W},
\]
with $\mathcal{Z}, \mathcal{W}$ supported in $Z \times Z$ and $W \times W$, respectively. We have the diagram:
\[
\begin{tikzcd}[column sep=small]
\ldots \arrow{r} & \text{H}^5_Z(X,\Z) \arrow{r}\arrow{d} & \text{H}^5(X,\Z) \arrow{r}{j^*_U}\arrow{d}{(\Gamma_{\phi})_* } & \text{H}^5(U,\Z) \arrow{r}\arrow{d}{\simeq}\arrow[dd,bend left=30,equal ] & \text{H}^6_Z(X,\Z)\arrow{d}{} \arrow{r} &\ldots \\
\ldots \arrow{r} & \text{H}^5_{W}(\mathbb{P}^5,\Z) \arrow{r} \arrow{d}{}&  \text{H}^5(\mathbb{P}^5,\Z) \arrow{r}{j^*_V} \arrow{d}{(\Gamma_{\phi^{-1}})_* } &  \text{H}^{5}(V,\Z) \arrow{r}\arrow{d}{\simeq} &  \text{H}^6_{W}(\mathbb{P}^5,\Z) \arrow{r}\arrow{d}{} & \ldots \\
\ldots \arrow{r} & \text{H}^5_Z(X,\Z) \arrow{r} & \text{H}^5(X,\Z)\arrow{r}{j^*_U}  & \text{H}^5(U,\Z) \arrow{r}& \text{H}^6_Z(X,\Z) \arrow{r} &\ldots \\
\end{tikzcd}
\vspace{-5mm}
\]
Since $\text{H}^5(\mathbb{P}^5,\Z) = 0$ then $\ker(j_U^*)=\text{H}^5(X,\Z)$, and hence $\operatorname{N}^2\text{H}^5(X,\Z) \simeq \text{H}^5(X,\Z)$.
\end{proof}

Thanks to the fact that $\text{H}^5(X,\Z)$ is of coniveau $2$, we conclude as in \cite[Corollaire 14 and Proposition 16]{FT} that $\operatorname{Griff}^3(X) = 0$. Thus, $\CH^3(X)_{\text{hom}} \simeq \CH^3(X)_{\text{alg}}$, and by \cite[Théorème 6.3]{Beau}, $\CH^3(X)_{\text{alg}} \simeq J^5(X)$, which is a principally polarized abelian variety. Since $X$ contains a plane $H \simeq \mathbb{P}^2 \hookrightarrow X$, the map $\CH^3(X) \to \text{H}^6(X,\Z)$ is surjective, and the resulting section gives a split exact sequence:
\[
0 \to J^5(X) \to \CH^3(X) \to \Z \to 0.
\]

\subsubsection{Chow groups of \texorpdfstring{$X_{(2,2,2)} \subset \mathbb{P}^{2m+4}$}{} when \texorpdfstring{$m \geq 3$}{}}

Now consider $X = X_{(2,2,2)} \subset \mathbb{P}^{2m+4}$ with $m > 2$. For $m \geq 2$, the variety $X$ is rational, with a birational map $\phi : X \dashrightarrow \mathbb{P}^{2m+1}$ induced by projection from a linear space $H$. Therefore, for all $i > 0$, the unramified cohomology groups $\text{H}^i_{\text{nr}}(X, \Q/\Z)$ vanish.

In general, for $i \neq m+1$, we have an isomorphism
\[
\CH^i(X) \simeq d_i \cdot \Z \oplus \CH^i(X)_{\text{tors}},
\]
with $d_i \in \N^{\geq 1}$ and $\CH^i(X)_{\text{hom}} \simeq \CH^i(X)_{\text{tors}}$. Since $X$ contains a hyperplane $H \subset X$ of dimension $m-1$, we have $d_i = 1$ for $m+2 \leq i \leq 2m+1$.

For $i = m+1$, we have a short exact sequence:
\[
0 \to \CH^{m+1}(X)_{\text{hom}} \to \CH^{m+1}(X) \to d_{m+1} \cdot \Z,
\]
which splits if and only if $d_{m+1} = 1$.

Note that $\CH^{2m}(X)$ is torsion-free: by \cite[Theorem 1]{BS83}, $\CH^{2m}(X)_{\text{hom}} = \CH^{2m}(X)_{\text{alg}} = 0$, hence $\CH^{2m}(X) \simeq \Z \cdot l$ where $l \simeq \mathbb{P}^1 \subset X$ is a line.

Following the proof of \cite[Corollaire 12]{FT}, one can observe that the diagonal of $X$ decomposes as
\[
\Delta_X = x \times X + l \times h + Z',
\]
where $x \in X$ is a point, $h = c_1(\mathcal{O}_X(1))$ is the hyperplane class, and $Z'$ is supported on $X \times T$ with $\operatorname{codim}_X(T) \geq 2$. Here, the constant $m$ in \cite[Corollaire 12]{FT} equals 1 due to the rationality of $X$ (see \cite[Proposition 6.1]{Par}).

\begin{corollary}\label{zero1}
Let $X = X_{(2,2,2)} \subset \mathbb{P}^{2m+4}$ be a smooth complete intersection with $m \geq 3$. Then for any integer $p \neq m$ and all $k \geq 1$, we have
\[
\operatorname{H}^p(X, \mathcal{H}^{p+k}(\Z)) = 0.
\]
\end{corollary}

\begin{proof}
The same argument used to prove Theorem~\ref{tors} shows that these groups are torsion.

Let $\phi : X \dashrightarrow \mathbb{P}^{2m+1}$ be the birational map, which is an isomorphism outside a linear subspace $H \subset X$ with $\operatorname{codim}_X(H) = m+2$. Then:
\[
\Delta_X = \phi^*(\Delta_{\mathbb{P}^{2m+1}}) + Z \in \CH^{2m+1}(X \times X),
\]
where $Z$ is supported on $H \times H$. Since $N_p = 1$ in the decomposition of the diagonal, the result follows.
\end{proof}

\begin{remark}
Since $X_{(2,2,2)} \subset \mathbb{P}^{2m+4}$ is rational, by \cite[Proposition 3.4]{CTV} we have for any commutative group $A$ and all $i \geq 0$ with $i \neq 2m+1$ that $\text{H}^i(X,\mathcal{H}^{2m+1}(A)) = 0$, and $\text{H}^{2m+1}(X,\mathcal{H}^{2m+1}(A)) \simeq A$.
\end{remark}

\begin{prop}
Let $X = X_{(2,2,2)} \subset \mathbb{P}^{2m+4}$ with $m \geq 3$. Then for all $0 \leq i \leq 2m+1$ with $i \neq m+1$, the groups $\CH^i(X)$ are torsion-free.
\end{prop}

\begin{proof}
The argument relies on the surjectivity of
\[
\text{H}^{i-1}_{\text{Zar}}(X,\mathcal{H}^i(\mu^{\otimes i}_{\ell^r})) \twoheadrightarrow \CH^i(X)[\ell^r].
\]
By Corollary~\ref{zero1}, for fixed $i \neq m+1$, we have $\text{H}^{i-1}_{\text{Zar}}(X,\mathcal{H}^i(\Z)) = 0$. Therefore, the group on the left vanishes, and hence $\CH^i(X)[\ell^r] = 0$ for all primes $\ell$ and all $r$, proving that $\CH^i(X)$ is torsion-free.
\end{proof}

\section{\texorpdfstring{\'Etale Chow groups of complete intersections with Hodge level one}{\'Etale Chow groups of complete intersections with Hodge level one}}\label{section:Chow etale level 1}

In this section we compute the \'etale Chow groups for smooth complete intersections in the projective space with Hodge structure of level one. More precisely, we will study the comparison maps
\[
\kappa^i: \operatorname{CH}^i(X)\to \operatorname{CH}_\et^i(X)
\]
that were introduced in \S \ref{section:Chow etale}. 

We will first discuss the case of Fano threefolds (which includes, in particular, cubic threefolds, quartic threefolds, and complete intersections of type $(2,3)$ in $\mathbb{P}^5$), and then turn to the case of higher-dimensional complete intersections.

\subsection{Fano threefolds}\label{subsec:Fano3folds}

Let $X$ be a smooth Fano threefold over the field of complex numbers $\C$. Such varieties satisfy the following two properties: Grabowski proved in \cite{Gra} that the integral Hodge conjecture holds for $X$, and that its Chow groups are representable. Using these facts, we can readily deduce the following isomorphisms:

\begin{theorem}
Let $X$ be a smooth Fano threefold. Then:
\[
\kappa^i:\CH^i_\et(X) \xrightarrow{\simeq} \CH^i(X) \textup{ for all }i\geq 0.
\]
\end{theorem}

\begin{proof}
For codimensions $0$ and $1$, the result holds for all smooth projective varieties over a field $k$ with $\operatorname{char}(k) = 0$. If $Y$ is a smooth projective variety over a field $k = \bar{k}$ of characteristic zero, then $\CH^{\dim(Y)}_\et(Y) \simeq \CH^{\dim(Y)}(Y)$. For the last étale Chow group, we consider the short exact sequence:
\[
0 \to \CH^2(X) \xrightarrow{\kappa^2} \CH^2_\et(X) \to \text{H}^3_{\text{nr}}(X, \Q/\Z(2)) \to 0.
\]
By \cite[Théorème 3.9]{CTV}, if $\CH_0(Y)$ is supported on a subvariety of $Y$ of dimension at most $2$, then the obstruction to the integral Hodge conjecture is isomorphic to $\text{H}^3_{\text{nr}}(X, \Q/\Z(2))$.

Furthermore, by \cite[Proposition 1.6]{Jann}, the condition that $\CH^2(X)_{\text{alg}}$ is representable is equivalent to the condition that $\CH_0(X)$ is supported on a subvariety of $X$ of dimension $1$. Since the integral Hodge conjecture holds for $X$, it follows that $\text{H}^3_{\text{nr}}(X, \Q/\Z(2)) = 0$, and therefore $\kappa^2$ is an isomorphism.
\end{proof}

\subsection{Cubic fivefolds}

Let $X\subseteq \mathbb{P}^6$ be a smooth cubic fivefold defined over $\C$. Its classical Chow groups are characterized as follows:
\begin{itemize}
\item $\CH^0(X)\simeq\Z$, generated by the fundamental class of $X$.
\item $\CH^1(X)\simeq\Z$, generated by the class of a hyperplane section $h:=c_1(\mathcal{O}_X(1))$.
\item $\CH^2(X)\simeq\Z$, generated by $h^2$. Since $\CH^2(X)_{\mathrm{hom}}=0$, the group $\CH^2(X)$ is isomorphic to the image of the cycle class map.
\item $\CH^3(X)\simeq\Z\oplus J^5(X)$, where $J^5(X)$ is a principally polarized abelian variety of dimension $21$ (see \cite[Théorème 1]{FT}).
\item $\CH^4(X)\simeq\Z$, generated by the class of a line contained in $X$.
\item $\CH^5(X)\simeq\Z$, generated by the class of a point, since $X$ is unirational.
\end{itemize}

\begin{theorem}
Let $X$ be a smooth cubic fivefold defined over $\C$. Then
\[
\kappa^i:\CH^i(X)\xrightarrow{\simeq}\CH^i_\et(X)\quad\text{for all }i\geq0.
\]
\end{theorem}

\begin{proof}
For $i=0,1,5$, the comparison maps are isomorphisms:
\[
\kappa^i:\CH^i(X)\xrightarrow{\simeq}\CH^i_\et(X).
\]
Since $\text{H}^3_{\mathrm{nr}}(X,\Q/\Z)=0$ (cf. \cite[Théorème 3]{FT}), it follows that $\kappa^2:\CH^2(X)\xrightarrow{\simeq}\CH^2_\et(X)$ is an isomorphism. For $i=4$, the comparison map $\kappa^4:\CH^4(X)\longrightarrow\CH^4_\et(X)$ is \emph{a priori} injective because $\CH^4(X)$ is torsion-free. The cokernel of $\kappa^4$ is isomorphic to
\[
\mathbb{H}^8_{\mathrm{Zar}}\!\left(X,\tau_{\geq6}R\pi_*\Z_X(4)_\et\right),
\]
and there exists a spectral sequence
\[
E^{p,q}_2:=\text{H}^p_{\mathrm{Zar}}\!\left(X,R^q\tau_{\geq6}R\pi_*\Z_X(4)_\et\right)\Longrightarrow\mathbb{H}^{p+q}_{\mathrm{Zar}}\!\left(X,\tau_{\geq6}R\pi_*\Z_X(4)_\et\right).
\]
Moreover,
\[
\text{H}^p_{\mathrm{Zar}}\!\left(X,R^q\tau_{\geq6}R\pi_*\Z_X(4)_\et\right)\simeq
\begin{cases}
\text{H}^p_{\mathrm{Zar}}\!\left(X,\mathcal{H}^{q-1}_\et(\Q/\Z(4))\right),&\text{if }q\geq6,\\
0,&\text{otherwise}.
\end{cases}
\]
From this and \cite[Théorème 3]{FT}, we deduce that the cokernel of $\kappa^4$ vanishes, hence $\kappa^4:\CH^4(X)\xrightarrow{\simeq}\CH^4_\et(X)\simeq\Z$ is an isomorphism. For $i=3$, the cokernel of $\kappa^3$ vanishes, so $\kappa^3$ is surjective. Its kernel is a subgroup of
\[
\mathbb{H}^5_{\mathrm{Zar}}\!\left(X,\tau_{\geq5}R\pi_*\Z_X(3)_\et\right),
\]
which, by the same reasoning, is trivial. Therefore, $\kappa^3:\CH^3(X)\xrightarrow{\simeq}\CH^3_\et(X)$ is an isomorphism.
\end{proof}

\subsection{Intersections of two and three quadrics}\label{subsec:quadrics-et}

In contrast to the classical setting, where one typically has explicit control of some of the Chow groups only, in the case of étale Chow groups we can compute \emph{all} of them. This is due to the fact that, for the varieties considered here, the rational Hodge conjecture holds, and we can apply the results of Rosenschon and Srinivas \cite{RS} to obtain a complete description in every codimension.

\begin{prop}
Let $X_{(2,2)}\subset\mathbb{P}^{2m+3}$ be a smooth intersection of two quadrics. Then the étale Chow groups of $X$ are given by:
\[
\CH_\et^i(X)\simeq
\begin{cases}
\Z[h^i], & \text{if } i\neq m+1,\\
\Z[h^{m+1}]\oplus J^{2m+1}(X), & \text{if } i=m+1,
\end{cases}
\]
where $h$ is a hyperplane section, and $J^{2m+1}(X)$ is the $(2m+1)$-intermediate Jacobian, which is isomorphic to the Jacobian of a smooth projective curve $C$.
\end{prop}

\begin{proof}
This follows directly from \cite{RS}. Since $\CH^i_\et(X)\to\mathrm{Hdg}^{2i}(X,\Z)$ is surjective if and only if $\CH^i(X)_{\Q}\to\mathrm{Hdg}^{2i}(X,\Q)$ is surjective, we have that $\CH^i_\et(X)\to\Z$ is surjective. Moreover, $\CH^i_\et(X)_{\mathrm{hom}}\simeq\CH^i_\et(X)_{\mathrm{tor}}$, and the torsion part of $\CH^i_\et(X)_{\mathrm{hom}}$ (by \cite[Proposition~5.1]{RS}) is isomorphic to $J^i(X)_{\mathrm{tors}}$. For the torsion-free part, we use that $\CH^i_\et(X)_\Q\simeq\CH^i(X)_\Q$ and that the torsion-free part of $J^{2m+1}(X)$ is uniquely divisible.
\end{proof}

\begin{prop}
Let $X_{(2,2,2)}\subset\mathbb{P}^{2m+4}$ be a smooth intersection of three quadrics. Then the étale Chow groups of $X$ are given by:
\[
\CH^i_\et(X)\simeq
\begin{cases}
\Z[h^i], & \text{if } i\neq m+1,\\
\Z[h^{m+1}]\oplus J^{2m+1}(X), & \text{if } i=m+1,
\end{cases}
\]
where $h$ is a hyperplane section, and $J^{2m+1}(X)$ is the $(2m+1)$-intermediate Jacobian of $X$.
\end{prop}

\begin{proof}
Again, this follows from \cite{RS}. The same argument as in the previous proposition applies: the surjectivity of $\CH^i_\et(X)\to\mathrm{Hdg}^{2i}(X,\Z)$ follows from its rational analogue, and the identification of $\CH^i_\et(X)_{\mathrm{hom}}$ with its torsion part reduces the problem to understanding $J^{2i-1}(X)_{\mathrm{tors}}$. The torsion-free part is obtained from the equality $\CH^i_\et(X)_\Q\simeq\CH^i(X)_\Q$.
\end{proof}

\begin{remark}
The above result, together with the discussion in \S \ref{subsec:Chow 2 quadrics} and \S \ref{subsec:Chow 3 quadrics}, implies that if $X$ is a $5$-dimensional intersection of two or three quadrics then $\kappa^i:\operatorname{CH}_{\et}^i(X)\xrightarrow{\simeq} \operatorname{CH}^i(X)$ are isomorphisms for all $i\geq 0$.
\end{remark}

\section{Decomposition of finite dimensional motives}\label{section:Motives}

In this section, we show that the previous results on the characterization of Chow groups (and \'{e}tale Chow groups) also hold in a motivic framework. The key point is the finite dimensionality of these motives in the sense of Kimura, see \cite{kim} and \cite[Chapters 4 and 5]{MNP} for details. This condition is more easily verified for varieties with Hodge structures of level~1. For instance, the motives of cubic threefolds and fivefolds are finite-dimensional \cite[Chapter~7, \S4]{huyB}, whereas the finite dimensionality of the cubic fourfold motive remains unknown.

More generally, we will observe that for the varieties that we are interested in, the motivic decomposition is given by the direct sum of Lefschetz motives and a $1$-motive given by the associated intermediate Jacobian, and thus the finite-dimensionality can be checked more easily.

We briefly recall, mostly in words, the category of Chow motives and some of its main properties, referring to \cite[Chapter 2]{MNP} for details. Let $\operatorname{SmProj}_k$ denote the category of smooth projective varieties over $k$. For $X,Y\in \operatorname{SmProj}_k$ with $X$ of pure dimension $d$, the correspondences of degree $r$ from $X$ to $Y$ are the elements of $\operatorname{Corr}^r(X,Y)_\Q:=\CH^{d+r}(X\times Y)_\Q$, which compose by the usual formula $g\circ f=(\operatorname{pr}_{XZ})_\ast(\operatorname{pr}_{XY}^\ast f\cdot \operatorname{pr}_{YZ}^\ast g)$. The category $\ch(k)_\Q$ of Chow motives has as objects the triples $(X,p,m)$, where $p\in \operatorname{Corr}^0(X,X)_\Q$ is an idempotent and $m\in\Z$, and the morphisms from $(X,p,m)$ to $(Y,q,n)$ are the correspondences of the form $q\circ \operatorname{Corr}^{n-m}(X,Y)_\Q\circ p$. The motive of $X$ is $h(X):=(X,\Delta_X,0)$, and a morphism $f\colon X\to Y$ induces $h(f)={}^t\Gamma_f\colon h(Y)\to h(X)$, where ${}^t\Gamma_f$ denotes the transpose of the graph of $f$.

\begin{remark}
    The category $\ch(k)_\Q$ is additive and pseudo-abelian (i.e.\ every idempotent morphism has both a kernel and a cokernel). Moreover, if we consider the triangulated category of motives $\operatorname{DM}(k,\Q)$, there is a fully faithful embedding of $1$-categories
    \[
        \ch(k)^{\mathrm{opp}}_\Q \hookrightarrow \operatorname{DM}(k,\Q).
    \]
    Another important point is that Chow motives appear in $\operatorname{DM}(k,\Q)$ as the weight-$0$ objects, with respect to the weight structure constructed by Bondarko \cite{Bon}.
\end{remark}

The category $\ch(k)_\Q$ is a tensor category, with $(X,p,m)\otimes(Y,q,n)=(X\times Y,p\otimes q,m+n)$, direct sums given by disjoint unions, and a duality $(X,p,m)^\vee=(X,{}^tp,d-m)$. We denote by $\mathbf{1}=(\operatorname{Spec}(k),\operatorname{id},0)$ the unit motive and by $\mathbf{L}=(\mathbb{P}^1_k,\mathbb{P}^1_k\times e,0)$ the Lefschetz motive, where $e$ is a rational point, so that $h(\mathbb{P}^1_k)=\mathbf{1}\oplus\mathbf{L}$ and $\mathbf{L}$ is inverse to the Tate motive $\mathbf{T}=(\operatorname{Spec}(k),\operatorname{id},1)$. The Chow groups and the cohomology, for a Weil cohomology theory $\operatorname{H}(-)$, of a motive $M=(X,p,m)$ are $\CH^i(M):=p_\ast\CH^{i+m}(X)_\Q$ and $\operatorname{H}^i(M):=p_\ast\operatorname{H}^{i+2m}(X)$. Finally, if $X$ has pure dimension $d$ and a rational point $e$, then $p_0(X)=e\times X$ and $p_{2d}(X)=X\times e$ are orthogonal projectors, and setting $h^0(X)=(X,p_0(X),0)$, $h^{2d}(X)=(X,p_{2d}(X),0)$, $p^+(X)=\Delta_X-p_0(X)-p_{2d}(X)$ and $h^+(X)=(X,p^+(X),0)$, one obtains a decomposition
\[
h(X)\simeq h^0(X)\oplus h^+(X)\oplus h^{2d}(X)\simeq \mathbf{1}\oplus h^+(X)\oplus\mathbf{L}^{\otimes d}.
\]

We now add a few remarks on the latter point. Fix a Weil cohomology theory $\operatorname{H}(-)$ and consider the diagonal class $\Delta_X$ of $X$. Let $\gamma_{X\times X}$ denote the cycle class map to $\operatorname{H}(-)$. By the Künneth formula one has
\[
    \gamma_{X\times X}(\Delta_X)=\sum_{i=0}^{2d}\Delta_i^{\mathrm{topo}}(X) \in H^{2d}(X\times X)\simeq \bigoplus_{i=0}^{2d} H^{2d-i}(X)\otimes H^i(X),
\]
with the convention that $\Delta_i^{\mathrm{topo}}(X) \in \operatorname{H}^{2d-i}(X)\otimes \operatorname{H}^i(X)$. Since $\operatorname{H}^{2d-i}(X)\otimes \operatorname{H}^i(X) \simeq \mathrm{End}(\operatorname{H}^i(X))$, we may regard $\Delta_i^{\mathrm{topo}}(X)$ as the identity on $\operatorname{H}^i(X)$. Thus, $\Delta_i^{\mathrm{topo}}(X)$ acts on $\operatorname{H}^*(X)$ as the identity on $\operatorname{H}^i(X)$ and as zero on all other summands of the above decomposition. This motivates the following definition.

\begin{defi}
    Let $X \in \mathrm{SmProj}_k$ be of pure dimension $d$. We say that $X$ admits a \emph{Chow-Künneth decomposition} if there exist projectors $p_i(X) \in \mathrm{Corr}^0(X,X)_\Q$ for $0\leq i \leq 2d$ such that:
    \begin{enumerate}
        \item $\sum_{i=0}^{2d}p_i(X) = \Delta_X$.
        \item The projectors $p_i(X)$ are pairwise orthogonal, i.e.
        \[
            p_i(X)\circ p_j(X)=
            \begin{cases}
                0 & \text{if } i\neq j,\\
                p_i(X) & \text{if } i=j,
            \end{cases}
        \]
        \item $\gamma_{X\times X}(p_i(X))=\Delta_i^{\mathrm{topo}}(X)$.
    \end{enumerate}
\end{defi}

It is conjectured that every smooth projective variety $X$ admits a Chow-Künneth decomposition. A result of Murre shows that, for a smooth projective variety $X \in \mathrm{SmProj}_k$, there exist projectors $p_1(X)$ and $p_{2d-1}(X)$ with the following properties (see also \cite[Chapter 6]{MNP}):

\begin{theorem}[{\cite{murre}}]
    Let $X \in \mathrm{SmProj}_k$ be smooth projective of dimension $d$. Then there exist $p_1(X), p_{2d-1}(X) \in \mathrm{Corr}^0(X,X)_\Q$ such that:
    \begin{enumerate}
        \item $p_1(X)$ and $p_{2d-1}(X)$ are orthogonal to each other and to $p_0(X)$ and $p_{2d}(X)$. They lift the Künneth components $\Delta_1^{\mathrm{topo}}$ and $\Delta_{2d-1}^{\mathrm{topo}}$, respectively, and satisfy $p_{2d-1}(X)={}^tp_{1}(X)$.
        \item Let $h^1(X)=(X,p_1(X),0) \in \operatorname{Chow}(k)_\Q$. Then
        \[
            \operatorname{H}^i(h^{1}(X))=
            \begin{cases}
                0 & \text{if } i\neq 1,\\
                \operatorname{H}^{1}(X) & \text{if } i=1,
            \end{cases}
            \qquad
            \CH^i(h^{1}(X))_\Q=
            \begin{cases}
                0 & \text{if } i\neq 1,\\
                \mathrm{Pic}^0(X)_\Q & \text{if } i=1.
            \end{cases}
        \]
        \item Let $h^{2d-1}(X)=(X,p_{2d-1}(X),0) \in \operatorname{Chow}(k)_\Q$. Then
       \begin{align*}
\operatorname{H}^i(h^{2d-1}(X)) &=
\begin{cases}
0 & \text{if } i\neq 2d-1,\\
\operatorname{H}^{2d-1}(X) & \text{if } i=2d-1,
\end{cases}\\[0.8em]
\CH^i(h^{2d-1}(X))_\Q &=
\begin{cases}
0 & \text{if } i\neq 2d-1,\\
\mathrm{Alb}(X)_\Q & \text{if } i=2d-1.
\end{cases}
\end{align*}
    \end{enumerate}
\end{theorem}

\begin{example}
    The Chow-Künneth decomposition is known to exist in the following cases\footnote{This list is not exhaustive.}: curves, surfaces (\cite{murre}), products of varieties admitting Chow-Künneth decompositions, abelian varieties, uniruled threefolds, complete intersections in projective space, Calabi-Yau threefolds (see \cite[Appendix~C]{MNP}), and double covers of $\mathbb{P}^n$ with $\dim \operatorname{H
    }^n(X,\Q)>1$ (\cite[Proposition~3.6]{LatM}), among others.
\end{example}

\begin{remark}
    Using the triangulated category of étale motives $\operatorname{DM}_\et(k,\Z)$ one can construct an étale analogue of the category $\ch(k)$, denoted $\ch_\et(k)$, which embeds fully faithfully into $\operatorname{DM}_\et(k,\Z)$. Its objects are triples $(X,p,m)$ with $X \in \operatorname{SmProj}_k$, $p \in \operatorname{Corr}^0_\et(X, X)$ satisfying $p \circ p=p$, and $m \in \Z$. Morphisms between two objects $(X,p,m)$ and $(Y,q,n)$ are given by
    \[
        \mathrm{Hom}_{\ch_\et(k)}\big((X,p,m),(Y,q,n)\big) = q \circ \operatorname{Corr}^{\,n-m}_\et(X,Y) \circ p.
    \]
     
    In this case we also have a contravariant functor $h_\et$:
    \begin{align*}
        h_\et: \operatorname{SmProj}_k &\to  \ch_\et(k) \\
        X   &\mapsto h_\et(X):=(X,\Delta_X,0) \\
        f:X \to Y &\mapsto h_\et(f)={}^t\Gamma_f:h_\et(Y)\to h_\et(X).
    \end{align*}
    
    For details on the construction of this category and the proof of the functoriality of étale correspondences, see \cite{RoSo}.    
\end{remark}

We now return to the decomposition of motives for varieties with Hodge structures of level~1, focusing in particular on varieties that are complete intersections. In all these cases, we have Chow-Künneth decomposition.

\begin{theorem}
Let $m \in \N$ and let $X$ be a smooth projective variety with $\operatorname{dim}_\C(X)=2m+1$ in one of the following cases:
\begin{enumerate}
    \item a complete intersection of two quadrics in $\Pro^{2m+3}$;
    \item a complete intersection of three quadrics in $\Pro^{2m+4}$ \emph{(}with smooth discriminant curve\emph{)};
    \item a cubic fivefold in $\Pro^6$;
    \item a Fano threefold.
\end{enumerate}
Then 
\[
    h_\et^{2m+1}(X)\simeq h^1_\et(J^{2m+1}(X))\otimes \mathbf{L}^m.
\]
\end{theorem}

\begin{proof}
The statement compiles several known results, each corresponding to one of the cases for $X$. All arguments and properties involved (notably finite dimensionality and motivic decompositions) can be found in the following references:
\begin{enumerate}
    \item For complete intersections of two quadrics, see \cite[Theorem 4.8]{reid} and \cite[Theorem 2.1]{lat}.
    \item For complete intersections of three quadrics, see \cite[Theorem 1]{NS} and \cite[Theorem 4.3]{boua}.
    \item For cubic fivefolds in $\Pro^6$, the result follows from the criterion in \cite[Theorem 4]{via}, together with the fact that $\CH_\et^i(X)_\Q \simeq \Q$ for all $0\leq i \leq 5$.
    \item For Fano threefolds, the statement follows from \cite[Theorems 8 and 15]{GG} concerning the representability of zero-cycles.
\end{enumerate}

If $X$ is one of the above cases, then the motive $h(X)\in \ch(\C)_\Q$ is finite dimensional in the sense of Kimura. By \cite[Corollary 5.6.4]{MNP} and the non-existence of phantom motives\footnote{i.e., non-trivial Chow motives which are homologically trivial.}, we obtain an integral étale Chow-Künneth decomposition of $h_\et(X)$ with
\[
    h_\et^{2n-1}(X):=(X,p), \qquad p=\Delta_X-\sum_{i=0}^{2n-2} p_{2i}^\et(X).
\]
This reflects the fact that there are no transcendental cohomology classes in codimensions other than $2m+1$. In each case, the Hodge decomposition
\[
    \operatorname{H}^{2m+1}_B(X,\C)=\bigoplus_{i+j=2m+1}\operatorname{H}^{i,j}(X)
\]
has nonzero components only in $\operatorname{H}^{m+1,m}(X)$ and $\operatorname{H}^{m,m+1}(X)$. Moreover, the intermediate Jacobians $J^{2m+1}(X)$ are abelian varieties. Thus, there exists a Chow-Künneth decomposition of the motive 
\[
    h^1_\et(J^{2m+1}(X))\otimes \mathbf{L}^m \to h_\et(X),
\]
so in particular we obtain a morphism $h^1_\et(J^{2m+1}(X))\otimes \mathbf{L}^m \to h_\et(X)$. The realization of the motive $h^1_\et(J^{2m+1}(X))\otimes \mathbf{L}^m$ corresponds to the first homology group of $J^{2m+1}(X)$. Indeed, there is an isomorphism
\[
    \operatorname{H}^B_1(J^{2m+1}(X),\Z)(-m)\simeq \operatorname{H}^{2m+1}_B(X,\Z).
\]

Explicitly, the $1$-motive associated to $J^{2m+1}(X)$ can be written as
\[
    h^1_\et(J^{2m+1}(X))\otimes \mathbf{L}^m=(J^{2m+1}(X),q,-m),
\]
where $q \in \operatorname{End}(h_\et(J^{2m+1}(X)))$ is the idempotent defining $h^1_\et(J^{2m+1}(X))\otimes \mathbf{L}^m$ as a direct factor of $h_\et(J^{2m+1}(X))\otimes \mathbf{L}^m$. At the level of homological motives there is an isomorphism
\[
    h^{2m+1}(X)_\Q\simeq h^1(J^{2m+1}(X))_\Q\otimes \mathbf{L}^m \in \ch_{\operatorname{hom}}(k)_\Q.
\]
Since these motives are finite dimensional in the sense of Kimura, this isomorphism lifts to the category of Chow motives:
\[
    h^{2m+1}(X)_\Q\simeq h^1(J^{2m+1}(X))_\Q\otimes \mathbf{L}^m \in \ch(k)_\Q.
\]
As this is induced by an algebraic cycle, its integral image is generated by an integral étale class, by \cite[Proposition 3.1.7]{RoSo}. Therefore, the isomorphism lifts further to an integral isomorphism of étale Chow motives $h^{2m+1}_\et(X)\simeq h^1_\et(J^{2m+1}(X)) \otimes \mathbf{L}^m \in \ch_\et(k)$.
\end{proof}

\begin{remark}
We remark that both $J^{2i+1}_{\alg}(X)$ and $J^{2i+1}_H(X)$ are abelian varieties, and that the cup product induces a bilinear Riemann form on $J^{2i+1}_H(X)$.
\end{remark}

From now on we restrict to varieties $X \in \operatorname{SmProj}_\C$ for which the standard Lefschetz conjecture $B(X)$ holds (see \cite[Chapter 3]{MNP} and \cite{Kle68,Kle94}). Let us recall the statement. Let $Y\subset X$ be a smooth hyperplane section. Consider its cohomology class $c_\Q^1([Y])$, and define the Lefschetz operator
\[
    L:\operatorname{H}^i_B(X,\Q)\to \operatorname{H}^{i+2}_B(X,\Q), \qquad 
    \alpha \mapsto \alpha \cup c_\Q^1([Y]).
\]
Its $r$-fold iterate is denoted $L^r$. By the hard Lefschetz theorem we have an isomorphism
\[
    L^j:\operatorname{H}^{d-j}_B(X,\Q)\xrightarrow{\simeq}\operatorname{H}^{d+j}_B(X,\Q),
\]
whose inverse we denote by $\Lambda^j$. Clearly, $\Lambda^r$ can be written as $\Lambda\circ \cdots \circ \Lambda$ ($r$ times).

\begin{conj}[Conjecture $B(X)$]
The topological correspondence $\Lambda$ is algebraic, i.e.\ there exists $Z \in \CH^{d-1}(X\times X)_\Q$ such that 
\[
    c_\Q^{d-1}(Z)\cup (-)=\Lambda.
\]   
\end{conj}

The conjecture $B(X)$ is known in the following cases (see e.g. \cite{Lie68,Kle94, Ara06}):
\begin{enumerate}
    \item Curves and surfaces (consequences of the Lefschetz $(1,1)$ theorem).
    \item Abelian varieties.
    \item Varieties $X$ for which $\operatorname{H}^\ast(X)\simeq \operatorname{CH}^\ast(X)$.
    \item Complete intersections in products of projective spaces.
    \item Grassmannians.
    \item Uniruled threefolds and unirational fourfolds.
\end{enumerate}

\begin{lemma}\label{bXcovering}
The conjecture $B(X)$ holds when $X \to \mathbb{P}^n$ is a smooth double cover ramified over a hypersurface $B\subset \mathbb{P}^n$ of degree $2d$, with $d \in \N^{\geq 1}$.
\end{lemma}

\begin{proof}
The argument is analogous to the proof that conjecture $B(X)$ holds when $X$ is an ample hypersurface in projective space (see, e.g., \cite[Proposition 4.3]{Kle94}). As explained in  
\cite[Proposition 1.4.7.(vii)]{Kle68}, the key point is that the rational version of the weak Lefschetz theorem allows one to construct an algebraic class representing $\Lambda_X$ from the corresponding class in the ambient space. Since $X$ can be realized as an ample hypersurface in the weighted projective space $\mathbb{P}(1^{n+1},d)$, the argument carries over verbatim by invoking the weak Lefschetz theorem for rational cohomology of toric varieties due to Batyrev and Cox \cite{BC94} (cf.\ the recent preprint \cite{Rau25}, where the pullback map is described at the level of integral cohomology).
\end{proof}

 For $i \in \N^{\geq 1}$ with $i\neq 1$ and $d=\operatorname{dim}(X)$, consider the following coniveau-like filtration
\[
    \operatorname{N}^i_\et \operatorname{H}_B^{2i+1}(X,\C) = 
    \sum_{\substack{C \in \operatorname{SmProj}_{\C} \\ \Gamma \in \CH_\et^{i+1}(C\times X)}}
    \operatorname{Im}\!\left\{\Gamma_* :\operatorname{H}^1(C,\C)\to \operatorname{H}^{2i+1}(X,\C)\right\},
\]
where $C$ runs over all smooth projective curves. This defines a sub-Hodge structure of $\operatorname{H}^{2i+1}(X,\C)$ of coniveau $i$. In addition, we consider the maximal sub-Hodge structure of $\operatorname{H}^{i+1,i}(X)\oplus \operatorname{H}^{i,i+1}(X)$, which we denote by $V$. Each of these Hodge structures has an associated intermediate Jacobian, denoted $J^{2i+1}_{\alg}(X)$ and $J^{2i+1}_H(X)$, respectively. In these cases we quotient by the lattices 
\[
    \operatorname{N}^i\operatorname{H}^{2i+1}(X,\Z):=\operatorname{N}^i\operatorname{H}^{2i+1}(X,\C)\cap \operatorname{H}^{2i+1}_B(X,\Z), 
    \qquad 
    V_\Z:=V \cap \operatorname{H}^{2i+1}_B(X,\Z),
\]
in the corresponding definitions. Both constructions yield abelian varieties, since their tangent spaces are contained in $\operatorname{H}^{i+1,i}(X)\oplus \operatorname{H}^{i,i+1}(X)$. Conjecturally, they coincide by the generalized Hodge conjecture in weight $2i+1$ and level $1$.

\begin{theorem}\label{theodecomp}
Let $X$ be a smooth projective variety of dimension $d$ verifying the standard Lefschetz conjecture $B(X)$. Then for all $2\leq i\leq d-1$, the motive $h_1(J_{\alg}^{i+1}(X))\otimes \mathbf{L}^i$ is a direct factor of the motive $h^+(X)$, where $h^+(X)$ is associated with the projector 
\[
    p^+(X):=\Delta-p_0(X)-p_1(X)-p_{2d-1}(X)-p_{2d}(X).
\]
Moreover, the generalized Hodge conjecture holds in weight $2i+1$ and level $1$ for $X$ if and only if $h_1(J_H^{i+1}(X))\otimes \mathbf{L}^i$ is a submotive of $h_\et(X)$. 
\end{theorem}

\begin{proof}
The fact that the motives $h_1(J_{\alg}^{i+1}(X))\otimes \mathbf{L}^i$ are subobjects of $h(X)\in \ch(k)_\Q$ follows from \cite{via}. We briefly recall the construction of the relevant projectors.

Consider $\operatorname{N}^i\operatorname{H}^{2i+1}(X,\C)$ and $V$ as the sub-Hodge structure of $\operatorname{H}^{2i+1}(X,\C)$ described above. By definition, there exists a curve $C=\amalg_j C_j$ (not necessarily connected) such that for some $\Gamma \in \CH^{i+1}_\et(C \times X)$ we obtain a surjection 
\[
    \Gamma_*: \operatorname{H}^1(C,\Z) \twoheadrightarrow \operatorname{N}^i\operatorname{H}^{2i+1}(X,\Z).
\]
Thus the induced map of abelian varieties $\Gamma_* : \operatorname{Alb}_C(k) \twoheadrightarrow J^{2i+1}_{\alg}(X)$ is also surjective. By the semisimplicity of the category of abelian varieties up to isogeny, there exists $\alpha \in \operatorname{Hom}_{\operatorname{AV}}(J^{2i+1}_{\alg}(X),\operatorname{Alb}_C(k))$ such that 
\[
    \Gamma_* \circ \alpha = m\cdot \operatorname{id}_{J^{2i+1}_{\alg}(X)}
\]
for some $m\in \Z$. Furthermore, we have a homomorphism with finite kernel $\Gamma^t_* \circ \Lambda: J^i_a(X)\hookrightarrow \operatorname{Pic}^0(C)$.

Now consider the dual morphism $\Gamma^t_*:J^{2i+1}_{\alg}(X)^\vee \to \operatorname{Pic}^0(C)$. For such a curve, let $\Theta$ be its theta divisor and define $\Gamma_{\lambda}:=\alpha \circ \lambda \circ \alpha^{\vee}$ for a polarization $\lambda:J^{2i+1}_{\alg}(X)^\vee  \to J^{2i+1}_{\alg}(X)$. For $\lambda$ such that $\Gamma_\lambda= \alpha \circ \Lambda^{-1}\circ \alpha^{\vee}$, set $\gamma:=\Gamma_\lambda$. Then define the projector
\[
    \pi_{2i+1,i}:= \Gamma_* \circ \gamma \circ \Gamma^t_* \circ L^{2i+1} \in \CH^d(X\times X)_\Q.
\]
By \cite[Theorem 2.8]{via} there is an isomorphism of motives 
\[
     \pi_{2i+1,i} \circ L^{2i+1} \circ \pi_{2(d-i)-1,d-i-1}: (X,\pi_{2i,i},0) \to (X,\pi_{2(d-i)-1,d-i-1},2i+1-d),
\]
while by \cite[Theorem 2.9]{via}, for such projectors $\pi_{2i+1,i+1}$ we obtain an isomorphism of motives 
\[
    (X,\pi_{2i+1,i+1})\simeq h_1(J^{2i+1}_{\alg}(X))(i).
\]
By the theory of \cite{RS}, these projectors and morphisms lift to $\ch_\et(k)$, so $h_1(J^{2i+1}_{\alg}(X))(i)$ is indeed a subobject of $h_\et(X)$. 

For the last statement, assume the generalized Hodge conjecture holds in level $1$ and weight $2i+1$. Then $J_{\alg}^{i+1}(X)\simeq J_H^{i+1}(X)$, hence $h_1(J^{2i+1}_H(X))(i)$ is a direct factor of $h_\et(X)$. Conversely, suppose that $h_1(J_H^{i+1}(X))(i)$ is a direct factor of $h(X)$. By the assumption on the Lefschetz conjecture, $h_1(J_H^{d-i}(X))(d-i)$ is also a direct factor of $h(X)$. Considering the injection of abelian varieties $i_J: J_{\alg}^{i+1}(X)\hookrightarrow J^{2i+1}_H(X)$, we obtain a map of motives 
\[
    i^*_J(i): h_1(J^{2i+1}_H(X))(i) \to h_1(J^{2i+1}_{\alg}(X))(i).
\]
Here $h_1(J_H^{i+1}(X))(i)=(J_H^{i+1}(X),p,i)$, where $p \in \operatorname{End}_{\ch(k)}(h(J_H^i(X)))$ is the associated projector. The realization of this $1$-motive gives
\[
    \operatorname{H}^*_B(h_1(J_H^{i+1}(X))(i),\Q) \simeq \operatorname{H}^1_B(J_H^{i+1}(X),\Q)(i) \simeq V_\Q^{\vee}(i),
\]
with both isomorphisms induced by algebraic cycles, due to the Lefschetz-type assumptions. Considering the injection $h_1(J_H^{i+1}(X))(i)\hookrightarrow h(X)$, we obtain an injection of Hodge classes $V_\Q\hookrightarrow \operatorname{H}^{i+1}_B(X,\Q)$. By semisimplicity, $V_\Q$ is a direct factor of $\operatorname{H}^{i+1}_B(X,\Q)$. By the hyperplane Lefschetz theorem there exists a curve $C\hookrightarrow J_H^{i+1}(X)$ such that 
\[
    h^i:\operatorname{H}_1^B(C,\Q(1))\xrightarrow{\simeq} \operatorname{H}_1^B(J_H^{i+1}, \Q(i)),
\]
yielding the commutative diagram
\[
\begin{tikzcd}
\operatorname{H}_1(C,\Q(1)) \arrow{dr}{i\circ h^i}[swap]{\simeq} \arrow{d} & \\
\operatorname{H}_1(J_H^{i+1}(X),\Q(i)) \arrow[r,hook] & \operatorname{H}^{i+1}_B(X,\Q(i)).
\end{tikzcd}
\]
This completes the proof.
\end{proof}

By Lemma~\ref{bXcovering}, Theorem~\ref{theodecomp} applies in particular to quartic double fivefolds, i.e.\ double covers of $\mathbb{P}^5$ branched along a smooth quartic hypersurface, which are studied in \cite{MRI}.


\begin{remark}\label{rmk:GM}
Among higher-dimensional Fano manifolds of Picard rank one, an important family whose Hodge theory and rationality have been extensively studied is given by Gushel--Mukai varieties \cite{DK18,DK20}. A Gushel--Mukai variety of dimension $n\in\{3,4,5,6\}$ is a smooth projective variety of the form
\[
X=\mathrm{CGr}(2,V_5)\cap \mathbb{P}(W_{n+5})\cap Q,
\]
where $V_5$ is a $5$-dimensional vector space, $\mathrm{CGr}(2,V_5)\subset \mathbb{P}(k\oplus\wedge^2V_5)$ is the cone over the Grassmannian $\mathrm{Gr}(2,V_5)$, $W_{n+5}\subset k\oplus\wedge^2V_5$ is a linear subspace of dimension $n+5$, and $Q\subset\mathbb{P}(k\oplus\wedge^2V_5)$ is a quadric hypersurface. In odd dimension, the middle cohomology of $X$ is a Hodge structure of level one \cite{DK20}. If $X$ is a Gushel--Mukai threefold or fivefold over $\C$, then all the comparison maps $\kappa^i\colon\CH^i(X)\to\CH^i_\et(X)$ are isomorphisms: for threefolds this follows from \S\ref{subsec:Fano3folds}, while for fivefolds the argument given for fivefold intersections of three quadrics in \S\ref{subsec:Chow 3 quadrics} and \S\ref{subsec:quadrics-et} applies, with the additional observation that $N=1$ in the decomposition of the diagonal of Theorem~\ref{Teovoi}.
\end{remark}

\printbibliography[title={Bibliography}]

\info
\end{document}